\pdfoutput=1

\documentclass[11pt, reqno, oneside, notitlepage]{article}

\usepackage[a4paper, total={6in, 9.1in}]{geometry}

\usepackage{amsmath,amscd}
\usepackage{amssymb}
\usepackage{amsthm}
\usepackage{comment}
\usepackage{graphicx}
\usepackage{epstopdf} 
\usepackage{mathrsfs}
\usepackage{cite}
\usepackage[ocgcolorlinks, linkcolor=blue]{hyperref}
\usepackage{authblk}

\usepackage{bm}
\usepackage{hyperref}
\usepackage{xcolor}
\hypersetup{
	colorlinks,
	linkcolor={blue},
	urlcolor={blue},
	citecolor={red}
}

\newtheorem{thm}{Theorem}[section]
\newtheorem{prop}[thm]{Proposition}

\newtheorem{lem}[thm]{Lemma}

\newtheorem{defi}[thm]{Definition}
\newtheorem{rmk}[thm]{Remark}

\newtheorem*{proposition*}{Proposition}

\numberwithin{equation}{section}

\title{On inverse problems in predator-prey models}

\author[1,*]{Yuhan Li}
\author[1,$\dagger$]{Hongyu Liu}
\author[1,$\natural$]{Catharine W. K. Lo}
\affil[1]{Department of Mathematics, City University of Hong Kong}
\affil[*]{yuhli2-c@my.cityu.edu.hk}
\affil[$\dagger$]{hongyu.liuip@gmail.com, hongyliu@cityu.edu.hk}
\affil[$\natural$]{wingkclo@cityu.edu.hk}
\date{}

\begin{document}
\maketitle

\begin{abstract}
In this paper, we consider the inverse problem of determining the coefficients of interaction terms within some Lotka-Volterra models, with support from boundary observation of its non-negative solutions. In the physical background, the solutions to the predator-prey model stand for the population densities for predator and prey and are non-negative, which is a critical challenge in our inverse problem study. We mainly focus on the unique identifiability issue and tackle it with the high-order variation method, a relatively new technique introduced by the second author and his collaborators. This method can ensure the positivity of solutions and has broader applicability in other physical models with non-negativity requirements. Our study improves this method by choosing a more general solution $(u_0,v_0)$ to expand around, achieving recovery for all interaction terms. By this means, we improve on the previous results and apply this to physical models to recover coefficients concerning compression, prey attack, crowding, carrying capacity, and many other interaction factors in the system.
Finally, we apply our results to study three specific cases: the hydra-effects model, the Holling-Tanner model and the classic Lotka-Volterra model.~
\end{abstract}

{\centering \section{INTRODUCTION}}
\subsection{Problem setup and background}
The predator-prey model is the building block of organisms and ecosystems as biomass grows from its resources. To struggle for their existence, species compete, evolve, and disperse. The model makes three simplified assumptions. First, the mutual relationship has only one predator and one prey. Second, if the number of predators drops below a certain threshold, the number of prey increases, while if the number of predators increases, the number of prey decreases. Third, if the number of prey reaches a certain threshold, the number of predators increases while the number of prey decreases. If the population of the prey is small, the number of predators decreases. According to their specific contexts, they can take forms such as resource-consumer, plant-herbivore, parasite-host, virus-immune system, and susceptible-infectious interactions. The relationships involved are general loss-win interactions; hence, we can generalize the models and apply them outside of ecosystems. In fact, when carefully studying seemingly competitive interactions, they are often disguised as some form of predator-prey interaction, which is how vital the predator-prey model is to us.

The theory of population dynamics in the prey-predator system was put up first by Lotka \cite{L1925}  and Volterra\cite{V1926} and has consequently been named the Lotka-Volterra (L-V) model. The ideas were soon  developed into qualitative applications and theoretical studies on biological mechanisms of the interrelations among populations. Due to the effort devoted by many scientists during these years, the system has a more precise form to adapt to different research objects. The well-posedness question, such as existence of the solution and its stability, is well-studied \cite{F980ec, DY1, Du97}. Furthermore, the L-V model helps us better understand those tipping points and critical transitions in ecosystems, which play a significant role in both the supervision of the change of the natural environment and the decision of an appropriate time to take human intervention. For example, to maintain the vegetation coverage and structural stability of plant communities in the pasture, governments advocate for rotational and intensive grazing methods, which avoids overgrazing and local overhunting. It is the calculations and predictions on thresholds that help researchers make plans for ranch management. One may see that the formation of spatial patterns caused by Turing instability has been used as a warning signal for dangerous critical points and impending critical transitions in complex ecosystems\cite{Reva}. We refer to Bazykin\cite{ B98} for the structural and dynamic stability of L-V models,  Turing\cite{T52} for his reaction-diffusion theory and Turing patterns, and Lu-Xiang-Huang-Wang\cite{LBi22} for their stability study on bifurcations in the diffusive Bazykin model.

Taking the intraspecific competition among the prey, the crowding effect for the predator, and the relationship between those two species into consideration, we have the following generalized L-V model\cite{F980ec}:
\begin{equation} \label{classic}
\begin{cases} 
\partial_{t}x=ax(1-\frac{x}{K})-yp(x,y),\\ 
\partial_{t}y=y(-c+dp(x,y))-hy^2, 
\end{cases} 
\end{equation} 
where $x(t),$ $y(t)$ stands for the population densities of the prey and predators at time $t$ respectively; $ax(1-\frac{x}{K})$ is the function reflecting the specific growth of the prey if there exists no predators, $a$ denotes the intrinsic growth rate, $K$ represents the carrying capacity; $h$ is a positive coefficient expressing the predator competition caused by self-limitation; $p(x,y)$ is the functional response depicting changes in the density of prey attacked by a predator per unit time.

There is a critical part of \eqref{classic}, which is the functional response $p(x,y).$ Functional responses have been studied extensively to describe different patterns for one predator responding to the changing density of its prey. According to Holling, functional responses are usually divided into three types, namely Holling Type I, Type II, and Type III.
In Type I there is a linear relation between the prey density and the maximum number of prey killed, an example is $p(x,y)=bx$ used by Bazykin in \cite{B76}. In Type II, the proportion of prey consumed declines monotonically with prey density, such as $p(x)=\frac{bx}{1+Ax}$ taken by Bazykin\cite{B76, B98}, Hainzl\cite{Ha92}, Lu-Huang\cite{LGl21}, and Kuznetsov\cite{Ku98}. Type III is described by an S-shaped relationship, where the proportion of prey consumed is positively correlated with some regions of prey density, see Freedman\cite{F97} with $p(x)=\frac{bx^2}{1+Ax^2}.$ The functional response helps evaluate two essential parameters: processing time (i.e., the time required for a predator to attack, consume, and digest prey) and attack rate or search efficiency (i.e., the speed at which predators search for prey). All of the above functional responses are modeled as a function of prey density only, which fails to model the competition among predators, no matter how large the predator density is. Instead, we also consider functional responses with both prey and predator dependence, such as ratio-dependent type functional response $p(x,y)=\frac{mx}{ax+by}$ contributed by Haque\cite{H09ratio} and Jiang et al.\cite{J21glo}, or the Beddington-DeAngelis type functional response $p(x,y)=\frac{mx}{ax+by+c}$ studied by Haque\cite{Haq11BD}.

Among the different forms of functional responses, Holling type II $p(x)=\frac{bx}{1+Ax}$ can give us the necessary information and can imitate a variety of well-known classic ecological modes. In the physical background, $b$ is the rate reflecting how fast the prey is killed by the predator and $A$ is the growing saturation effect of the predator when the population of the prey gets denser. Based on this functional response, we add the effects of predators and take two phenomena as examples: the hydra effect and the Holling-Tanner type.

Hydra effect occurs when the mean density of a species increases in response to greater mortality \cite{AbrHy05, AbrHy09}. Suppose it is predation that causes the increased mortality for a natural population. There are two steps in this interaction. Firstly, a species reduces its biomass as the mortality rate increases. Secondly, as a counterintuitive effect, the species enhance the equilibrium or time-averaged density. The interaction between predators and prey can be described by the following system of differential equations:
\begin{equation} \label{hydraeq}
      \begin{cases} 
      \partial_{t}u=au-bu-eu^2-(p+\lambda v)uv, \\ 
      \partial_{t}v=\mu (p+\lambda v)uv-mv
      \end{cases} 
    \end{equation} 
where $a$ and $b$ are the density-independent per capita birth and death rates of the prey respectively, $e$ is the intraspecific competition rate coefficient, $\mu$ is the efficiency of food conversion into offspring, $m$ is the per capita mortality rate of the predator, $p$ stands for the attack rate of an individual predator, and $\lambda$ is the strength of cooperation during hunting.

The Holling-Tanner model describes the dynamics of a generalist predator which feeds on its favourite food item as long as it is in abundant supply and grows logistically with an intrinsic growth rate and a carrying capacity proportional to the size of the prey. The system with Holling Tanner functional response is given by
\begin{equation} \label{clacha}
\begin{cases} 
\partial_{t}u(x,t)-d_1\Delta u(x,t)=u(1-u-bv/(1+mu)), &\  \text{in}\  Q,\\ 
\partial_{t}v(x,t)-d_2\Delta v(x,t)=v(d-v+cu/(1+mu)), &\  \text{in}\  Q,\\ 
\partial_{\nu}u(x,t)=\partial_{\nu}v(x,t)=0 &\  \text{on}\   \Sigma,\\ 
\end{cases} 
\end{equation} 
where all the coefficients are positive except $d.$ If $d>0,$ the predators can survive without the prey; otherwise, the predators will become extinct in the absence of the prey. Here we consider $\Omega$ as a bounded domain in $\mathbb{R}^N,$ $N \geq1,$ with smooth boundary $\partial \Omega.$ $\nu(x)$ is the unit outward normal vector on the boundary and $\Delta$ is the Laplace operator signifying the diffusion process. Denote $Q:=\Omega \times (0,\infty),$ $\Sigma:=\partial \Omega \times (0,\infty).$ This system has been shown to possess positive steady solutions by Du-Lou in \cite{DY1}. Moreover, if $m=0,$ Leung\cite{L78} proved that all the positive solutions of \eqref{clacha} converge to a constant steady state without the influence of initial data when time goes to infinity. A similar conclusion holds for small positive $m.$

Multiple results have been obtained for this L-V model discussing the stability of a constant steady state, the existence of non-constant positive steady states, and the Turing bifurcation with its pattern. Moreover, to reveal the nonlinear dynamics and complex bifurcation phenomena of \eqref{clacha}, Lu et al. in \cite{LBi22} put up some rigorous answers with the following system:
\begin{equation} \label{Difex}
\begin{cases} 
\partial_{t}u-D_1\Delta u= au(1-\frac{u}{K})-\frac{buv}{1+Au} &\  \text{in} \  Q,\\ 
\partial_{t}v-D_2\Delta v= -cv+\frac{duv}{1+Au}-hv^2 &\  \text{in}\   Q,\\ 
\partial_{\nu}u(x,t)=\partial_{\nu}v(x,t)=0 &\  \text{on}\   \Sigma,\\ 
u(x,0)= u_0(x) , \; v (x,0)= v_0(x) &\  \text{in}\   \Omega 
\end{cases} 
\end{equation} 
where $u(t),$ $v(t)$ stands for the population densities of the prey and predators at time $t$ respectively, corresponding to \eqref{classic}. $D_1, D_2$ are diffusion coefficients of the species,  $a$ depicts the reproduction rate of the prey population when there is no predator, $c$ is the natural mortality rate for the predator and $h$ is the coefficient reflecting the crowding effect. Then it is easy to see that $\frac{1}{A}$ is the prey density when the consumption of the predator is at half of its maximum number, and $\frac{b}{A}$ reports the maximum consumption of the predator. $\frac{d}{b}$ is the conversion efficiency among two species. All the parameters in \eqref{Difex} are positive. In later discussion, we shall recover these interaction coefficients from measurements of the boundary data of the population densities $u$ and  $v,$ and we need the positivity of the coefficients to ensure the positivity of $u$ and $v.$ 

In this paper, we consider the following general coupled semi-linear system of the predator-prey model:
\begin{equation}\label{moduse}
\begin{cases} 
\partial_{t}u-d_1\Delta u= F(x,t,u,v) &\  \text{in} \   Q,\\ 
\partial_{t}v-d_2\Delta v= G(x,t,u,v) &\  \text{in}  \  Q,\\ 
\partial_{\nu}u=\partial_{\nu}v=0&\  \text{on}\   \Sigma,\\ 
u(x,0)= f(x) , \; v (x,0)= g(x)&\  \text{in}\   \Omega.
\end{cases} 
\end{equation} 
Here, $F(x,t,p,q)$ and $G(x,t,p,q):$ $\Omega \times (0,T) \times \mathbb{R}^{n}\times \mathbb{R}^{n} \to \mathbb{R}$ are real-valued functions with respect to $p$ and $q.$ Both $u$ and $v$ are required to be non-negative to ensure their physical meaning. More details for $F(x,t,u,v)$ and $G(x,t,u,v)$ will be given in subsection \ref{Adc}.
 
Compared with direct studies on predator-prey models, the theory of this corresponding inverse problem is an emerging research area. Several numerical studies have been done on the inverse problem of biology, see examples as \cite{PW83, EW9}, but there are hardly any theoretical works on the inverse problems for the predator-prey model. Liu-Lo\cite{LL23} made a similar work for the inverse problem in parabolic equations and then applied their conclusion to biological models as examples. Ding-Liu-Zheng \cite{DL23} solved the unique identifiability issue for a fixed-form predator-prey model using monotonicity properties and comparison principle, which is different from our technique. Therefore, it is meaningful for us to apply high-order variation method to achieve unique identifiability results with population densities for both species being non-negative. Our results adapt for more complex forms for the L-V models, have fewer restrictions for preconditions, and can recover more coefficients.

In this paper, we mainly consider the inverse problem of determining the interaction coefficients, $F$ and $G,$ with $u(x,t)$ and $v(x,t)$ kept non-negative. 
To this end, we introduce a measurement map $\mathcal{M}^{+}_{F,G}$ as below:
\begin{equation}\label{map}
\mathcal{M}^{+}_{F,G}(f(x),g(x))=\big( (u(x,t),v(x,t))\vert_{\Sigma},u(x,T),v(x,T)\big),\    x \in \Omega,   \end{equation} 
where the sign $'+'$ stresses the non-negativity of the solutions to the coupled parabolic system. The measurement map tells that for a given pair of $F$ and $G,$ $\mathcal{M}^{+}_{F,G}$ maps from the initial population distribution $(f,g)$ to functions $u(x,t),$ $v(x,t)$ on the boundary. The details of $\mathcal{M}^{+}_{F,G}$ will be stated in section \ref{forp}. And we can formulate our inverse problem into:
\begin{equation}\label{ip}
\mathcal{M}^{+}_{F,G} \to F\, \, \text{and}\, \, G.\     \    \end{equation}

This inverse problem is realistically meaningful. As we know in \eqref{hydraeq} and \eqref{Difex}, each coefficient in the L-V model has its physical sense, which plays a role in human beings' supervision of the environment or offers guidance on how to facilitate a better change in ecology. Meanwhile, it is easier for professionals to detect population densities for species rather than calculating interactions among different creatures. Hence, we are motivated to explore the inverse problem \eqref{ip}. Specifically, we study the unique identifiability issue, which is a primary subject for an inverse problem. In mathematical words, it asks whether we can build such a one-to-one relationship between two measurement maps:
\begin{equation}\label{ideni}
\mathcal{M}^{+}_{F_1,G_1}=\mathcal{M}^{+}_{F_2,G_2} \  \   \text{if and only if} \  \   (F_1,G_1)=(F_2,G_2), \end{equation}
where $(F_j,G_j),$ $j=1,2,$ are two configurations. More details for this aspect will be given in section \ref{prelims}. 

\subsection{Technical development}
Though there is little research on inverse problems for the L-V model, we can find many inverse problems working on nonlinear parabolic equations, which are used to depict physical systems, see for example \cite{SR} for a work on a single nonlinear parabolic equation. Similar results have also been considered for systems, such as in \cite{CSE22,BAC09, PV}. There have also been some works studying inverse problems in biology, see \cite{DL23, LL23}. However, the results are very limited. In this paper, we improve on these previous results and can essentially achieve the recovery of all coefficients in a Bazykin prey-predator model. As such, we have a more comprehensive range of applications, such as the allowance for free-area research observation.

Two main challenges for this paper are the non-negativity constraint and homogeneous Neumann boundary condition. First, we focus on addressing the non-negativity constraints on solutions. Many works on the forward problem for parabolic systems, such as \cite{FA8, PV, SL20}, consider this issue. However,  articles enforcing non-negativity for the solutions in the inverse study are minimal, see \cite{LL23, LZ23, DL23}. Yet, since these models represent the population of the predator and prey species, the conclusions drawn from the theoretical study can only be physically realistic if the solutions are non-negative.

To address the issue of the positivity of the solutions, we primarily make use of the high-order variation method introduced in \cite{LL23,LZ23}, which is based on the method of successive linearization essential for treating nonlinear equations. The main idea used in \cite{LL23,LZ23} is to consider a specific form of linearization around zero-value for the input initial data. In those works, the authors set
\begin{equation}\label{characv}
u(x;\varepsilon)=\sum\limits^{\infty}_{l=1}\varepsilon^{l}f_{l} \  \  \text{on}  \,  \   \Omega \  \  \text{for} \,   \  f_1 \geq 0.
\end{equation}
Since $\varepsilon$ is a small positive variable, the positivity of $f_i(x),$ $i \geq 2$ would not impact the non-negativity of $u(x; \varepsilon)$ on the boundary.  Moreover, \cite{LZ23}  chose the special form $(0,1)$ and \cite{LL23} linearized around $(0,0).$  In this paper, we consider the case of a general initial input data $(u_0,v_0).$ Here, either one of the initial data or both data can cover the value of $0$ if the form of PDEs allows, which greatly broadens the usage for those equations with different known solutions. Compared to \cite{LL23}, which focuses on the recovery of Taylor coefficients for $G(x,t,u,v)$ about partial derivatives to $u,$ we let $G_u=0$ and recover all the rest of the coefficients. Meanwhile, we only have assumptions on the first-order Taylor coefficients of $F(x,t,u,v),$ which is less restrictive than \cite{LL23}. The forms of the source functions $F,G$ used in this paper also cover the form of $F$ in \cite{DL23}. Hence, our results recover more coefficients under comparatively loose requirements for more complex systems. It is one of the first few that uses the high-order variation method (only previously used in \cite{LL23,LZ23, DL23} to the best of our knowledge), and is currently the most versatile. As many real-life physical models have the non-negativity requirement, this technical development is widely applicable in tackling inverse problems.

One thing worth highlighting is that we show classic successive linearization comprehensively to manifest the advantages of introducing the high-order variation method into proofs. Except for the guarantee of non-negativity for the population densities, it also helps to simplify the calculation within the proofs, see sections \ref{highlinear}, \ref{highvary}. 

Constructing Complex-Geometric-Optics (CGO) solutions is customarily applied to studying the inverse problem with coupled parabolic equations. However, there are two difficulties for us in using this method to solve \eqref{moduse}:  the requirement of satisfying Neumann boundary conditions and the challenge of building non-negative CGO solutions. Our proposed method precisely avoids constructing such CGO solutions, and the core idea is to separate the time and space variables by adopting the method used in \cite{LZ23}. \cite{DL23} also offers another technique by a different selection of admissible classes.

The remainder of this paper is organized as below. We provide the basic notations and present the main result in section \ref{prelims}, and the well-posedness of the forward problem is covered in section \ref{forp}. Detailed descriptions of the linearization methods are given in section \ref{highva}, along with some theoretical preconditions for the proof of our main theorem \ref{mainthm}, which is verified in section \ref{proofm}. We apply our results to biological models in section \ref{appbio}.\\

{\centering \section{PRELIMINARIES AND STATEMENT OF MAIN RESULTS}  \label{prelims} }

\subsection{Notations and basic settings}
First, we introduce the norm of $C^{k+\alpha}(\bar{\Omega}).$ $C^{k+\alpha}(\bar{\Omega})$ is the subspace of $C^{k}(\bar{\Omega})$ such that for $k \in \mathbb{N}$ and $0<\alpha <1,$ we say that a function $\phi \in C^{k+\alpha}(\bar{\Omega})$ if for any $|l| \leq k,$ $D^{l} \phi$ exists and are H\"older continuous with exponent $\alpha.$ The norm is defined as
\[ \Vert \phi \Vert_{C^{k+\alpha}(\bar{\Omega})} =\sum\limits_{\vert l \vert \leq k} \Vert D^l \phi \Vert_{\infty}+\sum\limits_{\vert l \vert=k} \mathop{\text{sup}}\limits_{x \neq y}\frac{\vert D^l \phi(x)-D^l \phi(y)\vert}{\vert x-y \vert ^{\alpha}}.\]

For functions depending on both the time and space variables,  we define its space as $C^{k+\alpha, \frac{k+\alpha}{2}}(Q).$ We say that a function $\phi \in  C^{k+\alpha, \frac{k+\alpha}{2}}(Q)$ if $D^{l}D^{j}_{t} \phi$ exists and are H\"older continuous with exponent $\alpha$ in the space variable and exponent $\frac{k+\alpha}{2}$ in the time variable for all $l \in \mathbb{N}^n,$ $j \in \mathbb{N},$ $\vert l \vert + 2j \leq k.$ The norm is defined as
\[ \begin{split} 
\Vert \phi \Vert_{C^{k+\alpha, \frac{k+\alpha}{2}}(Q)} := & \sum\limits_{\vert l \vert + 2j \leq k} \Vert D^{l}D^{j}_{t}\Vert_{\infty}+\sum\limits_{\vert l \vert +2j=k} \mathop{\text{sup}}\limits_{t,x \neq y}\frac{\vert  \phi(x,t)- \phi(y,t)\vert}{\vert x-y \vert ^{\alpha}}+\\
&\sum\limits_{\vert l \vert+ 2j=k} \mathop{\text{sup}}\limits_{t \neq t^{'}, x}\frac{\vert  \phi(x,t)- \phi(x,t^{'})\vert}{\vert t-t^{'} \vert ^{\alpha/2}}.
\end{split} \]

The spaces $H^s(O)$ and $L^2(0,T;H^s(O)),$ $O=\Omega, \Sigma$ are standard Hilbert and Bochner spaces for  $s \in \mathbb{R}$ respectively. We introduce the following function space,
\[H_{\pm}(Q):=\{u\in \mathcal{P}^{'}(Q) \vert u \in L^2(Q)\, \text{and} \, (\pm \partial_{t}-d_1\Delta +q)u \in L^2(Q) \},\]
with its norm defined using its graph norm as 
\[\Vert u \Vert^2_{H_{\pm}(Q)}:=\Vert u \Vert^2_{L^2(Q)}+\Vert (\pm \partial_{t}-d_1\Delta +q)u \Vert^2_{L^2(Q)}.\]

Next, we restrict our study of the nonlinear system of \eqref{moduse} to the following form:
\begin{equation}\label{modac}
\begin{cases} 
\partial_{t}u-d_1\Delta u= F(x,t,u,v) &\  \text{in} \   Q,\\ 
\partial_{t}v-d_2\Delta v= G(x,t,u,v) &\  \text{in}  \  Q,\\ 
\frac{\partial u}{\partial \nu}=\frac{\partial v}{\partial \nu}=0&\  \text{on}\   \Sigma,\\ 
u(x,0)= f(x) \geq 0 , \; v (x,0)= g(x) \geq 0 &\  \text{in}\   \Omega ,\\
u(x,t) \geq 0 , \; v (x,t) \geq 0 &\  \text{in}\   Q.
\end{cases}
\end{equation} 
The functions $F(x,t,p,q), G(x,t,p,q)$ are analytic with respect to $p$ and $q,$ and are of the forms below:
\begin{equation}\label{Fr}
   F(x,t,p,q):=\sum\limits^{\infty}_{\substack{m\geq1,\,  n\geq 0,\, m+n\geq 2,\\ h=0,1}}\alpha_{mnh}\frac{p^m q^n}{1+p^h}, 
\end{equation}
and
\begin{equation}\label{Gr}
   G(x,t,p,q):=\sum\limits^{\infty}_{\substack{m\geq0,\,  n\geq 1, \,m+n\geq 2,\\ h=0,1} }\beta_{mnh}\frac{p^m q^n}{1+p^h}.
\end{equation}
In particular, this means that \eqref{modac} includes the physical models \eqref{clacha} and \eqref{Difex} that we discussed previously. The Neumann conditions given here illustrate that the domain is an enclosed habitat with reflecting boundary, which causes a great technical difficulty for us in building CGO solutions for \eqref{modac}, as discussed in the previous section. 

\subsection{Admissible class}\label{Adc}
Suppose $(u_0,v_0)$ is a known non-negative constant solution of \eqref{modac} and $F,G$ are analytic. Then we can introduce the admissible classes of the source functions in the L-V model.

\begin{defi}\label{admissible1}
We say that $U(x,t,p,q): \mathbb{R}^n \times \mathbb{R} \times \mathbb{C}  \times \mathbb{C} \to \mathbb{C}$ is admissible, denoted by $U \in \mathcal{A},$ if:

(a) The map $(p,q) \mapsto U(\cdot, \cdot, p,q) $ is holomorphic with value in $C^{2+\alpha, 1+\frac{\alpha}{2}}(\bar{Q}),$

(b) $U(x,t,u_0,v_0)=0$ for all $(x,t) \in Q,$ 

(c) $U^{(0,1)}(x,t,\cdot,v_0)=U^{(1,0)}(x,t,u_0,\cdot)=0$ for all $(x,t) \in Q,$ 

(d) Taylor coefficients of all orders for $U$ are constants.

It is clear that if $U$ satisfies these four conditions, it can be expanded into a power series
\[U(x,z,p,q)=\sum\limits^{\infty}_{m,n=1}U_{mn}\frac{p^m q^n}{(m+n)!},\]
where $U_{mn}=\frac{\partial^m}{\partial p^m}\frac{\partial^n}{\partial q^n}U(x,t,u_0,v_0)$ is a constant.
\end{defi}

\begin{defi}\label{admissible2}
We say that $V(x,t,p,q): \mathbb{R}^n \times \mathbb{R} \times \mathbb{C}  \times \mathbb{C} \to \mathbb{C}$ is admissible, denoted by $V \in \mathcal{B},$ if: 

(a) The map $(p,q) \mapsto V(\cdot, \cdot, p,q) $ is holomorphic with value in $C^{2+\alpha, 1+\frac{\alpha}{2}}(\bar{Q}),$

(b) $V(x,t,u_0,v_0)=0$ for all $(x,t) \in Q,$ 

(c) $V^{(1,0)}(x,t,u_0,\cdot)=0$ for all $(x,t) \in Q,$ 

(d) Taylor coefficients of all orders for $V$ are constants.

It is clear that if $V$ satisfies these four conditions, it can be expanded into a power series
\[V(x,z,p,q)=\sum\limits^{\infty}_{m\geq1,\, n\geq 0}V_{mn}\frac{p^m q^n}{(m+n)!},\]
where $V_{mn}=\frac{\partial^m}{\partial p^m}\frac{\partial^n}{\partial q^n}V(x,t,u_0,v_0) $ is a constant.
\end{defi}

It can be easily seen that for $F\in \mathcal{A}$ and $G\in \mathcal{B},$ $F$ and $G$ are of the forms \eqref{Fr} and \eqref{Gr} respectively.

\begin{rmk}\label{illuad}
\textit{By extending functions} $F$ \textit{and} $G$ \textit{of real variables to the complex plane of} $z$-\textit{variables, the above admissible conditions have a prior effect on these functions. Specifically, they are given by} $ \tilde{F} ( \cdot, \cdot, p, q) $ \textit{and} $ \tilde{G} (\cdot,\cdot, p, q)$\textit{ respectively, and their functions as complex variables }$p, q$ \textit{are holomorphic. Therefore,} $F$ \textit{and} $G$ \textit{are the restrictions of} $\tilde{F}$ \textit{and} $\tilde{G}$ \textit{ to the real line.} 
\end{rmk}

\begin{rmk}\label{rmk24}
The definitions for admissible classes $\mathcal{A}$ and $\mathcal{B}$ help to show great difference for us from \cite{LL23}. Though we choose not to recover $G_u,$ we can recover more coefficients in higher-order variation system. We also simplify the non-negative restriction by this setting, see section \ref{proofm}.
 \end{rmk}
 
\subsection{Main uniqueness identifiability results}
In this part, we propose our major result for the inverse problems, which shows in a generic scenario that we can uniquely recover the coefficients functions $F$ and $G$ from the measurement map $\mathcal{M}^{+}_{F,G}. $ We state our conclusions into the following theorem.
\begin{thm} \label{mainthm}
Assume $F_j \in \mathcal{A},$ $G_j \in \mathcal{B},$  $ (j=1,2). $ Let $\mathcal{M}^{+}_{F_j, G_j}$ be the measurement map associated to the following system:
\begin{equation}
\begin{cases} 
\partial_{t} u_{j}-d_1\Delta u_j=F_j(x,t,u,v), &\  \text{in}\  Q,\\ 
\partial_{t}  v_{j}-d_2\Delta v_j= G_j(x,t,u,v), &\  \text{in}\  Q,\\ 
\partial_{\nu} u_{j}(x,t)=\partial_{\nu} v_{j}(x,t)=0, &\  \text{on}\  \Sigma,\\ 
u_j(x,0)= f(x) , \; v_j (x,0)= g(x),& \  \text{in} \ \Omega. 
\end{cases} 
\end{equation} 
If for any $f,g \in C^{2+\alpha}(\Omega),$ one has  \[\mathcal{M}^{+}_{F_1,G_1}(f,g)=\mathcal{M}^{+}_{F_2,G_2}(f,g),\]
then it holds that 
\[F_1(x,t,u,v)=F_2(x,t,u,v),\  \, G_1(x,t,u,v)=G_2(x,t,u,v)\  \ \text{in}\  \, Q.\]
\end{thm}

Physically speaking, our results can be explained as we input a certain population of prey from the reflecting boundary into a bounded region in $\Omega,$ and measure the initial population densities of predator and prey. If all these values are the same, then we can recover the interaction coefficients, i.e. the maximum consumption of the predator, and the conversion efficiency between two species and prey density when the consumption of the predator is at half of its maximum quantity.~\\

{\centering \section{WELL-POSEDNESS OF THE FORWARD PROBLEMS}   \label{forp}}
This section aims to study the local and global well-posedness for initial-boundary value problems of semilinear parabolic equations. We consider the following equation:
\begin{equation}\label{we1}
\begin{cases} 
\partial_{t}u-d_1\Delta u+b_1(x,t,u,v)= 0, &\  \text{in}\  Q,\\ 
\partial_{t}v-d_2\Delta v+b_2(x,t,u,v)= 0, &\  \text{in}\  Q,\\ 
\partial_{\nu}u=\partial_{\nu}v=0, &\  \text{on}\  \Sigma,\\ 
u(x,0)=\tilde{f}(x), &\  \text{in}\  \Omega,\\
v(x,0)=\tilde{g}(x), &\  \text{in}\  \Omega,\\  
\end{cases} 
\end{equation} 
where $d_1,d_2$ are positive constants,  $\tilde{f}, \tilde{g}\in C^{2+\alpha}(\bar{\Omega}).$ Moreover, $b_j(x,t,u,v),$ $j=1,2,$ should also satisfy the following condition:
\begin{equation}\label{we1b}
b_j \in C^2(\bar{Q}\times \mathbb{R}) \  \text{and} \   b_j(\cdot,\cdot,u_0,v_0)=0,\  b_j(x,t,\cdot, \cdot)=0 \   \text{in}  \   Q,\  j=1,2.  \end{equation} 

Based on \cite{LSU88}, we recall the well-posedness result and Schauder estimates for linear parabolic equations. 

\begin{lem}\label{lemma31}
For $\alpha \in (0,1),$ assume that $p,q \in C^{2+\alpha}(\bar{Q}),$ $h_1,h_2 \in C^{\alpha,\alpha/2}(\bar{Q}).$ For any $\tilde{f},\tilde{g} \in C^{2+\alpha}(\bar{\Omega})$ with the compatibility conditions:
\begin{equation}\tilde{f}(x)=0\   \text{and}\  d_1 \Delta \tilde{f}(x)-p(x,0)\tilde{f}(x)+h_1(x,0)=0 \  \text{on}\  \Gamma,\end{equation}
\begin{equation}\tilde{g}(x)=0\   \text{and}\  d_2 \Delta \tilde{g}(x)-q(x,0)\tilde{g}(x)+h_2(x,0)=0 \  \text{on}\  \Gamma.\end{equation}
The following linear parabolic system:
\begin{equation}
\begin{cases} 
\partial_{t}u-d_1\Delta u+pu= h_1, &\  \text{in}\  Q,\\ 
\partial_{t}v-d_2\Delta v+qv= h_2, &\  \text{in}\  Q,\\ 
\partial_{\nu}u=\partial_{\nu}v=0, &\  \text{on} \  \Sigma,\\ 
u(x,0)=\tilde{f}(x), &\  \text{in}\  \Omega,\\
v(x,0)=\tilde{g}(x), &\  \text{in}\   \Omega,\\  
\end{cases} 
\end{equation} 
admits a unique solution $(u,v)\in C^{2+\alpha,1+\alpha/2}(\bar{Q}).$ And we obtain the below estimate:
\[ \begin{split}
\Vert u\Vert_{C^{2+\alpha,1+\alpha/2}(\bar{Q})}+\Vert v\Vert_{C^{2+\alpha,1+\alpha/2}(\bar{Q})} \leq C( & \Vert \tilde{f}\Vert_{C^{2+\alpha}(\bar{Q})}+\Vert \tilde{g}\Vert_{C^{2+\alpha}(\bar{Q})}+\\
&\Vert h_1\Vert_{C^{\alpha,\alpha/2}(\bar{Q})}+\Vert h_2\Vert_{C^{\alpha,\alpha/2}(\bar{Q})}).
\end{split} \]
\end{lem}

It is apparent to see that if $h_1=h_2=0$ in $Q,$ $\tilde{f}=\tilde{f}_{x_i}=\tilde{f}_{x_i x_j}=0,$ $\tilde{g}=\tilde{g}_{x_i}=\tilde{g}_{x_i x_j}=0,$ $(i,j=1,\cdots,n)$ on $\Gamma,$ then the compatibility condition \eqref{we1b}  holds true. 

Moreover, we can obtain the following local well-posedness for system \eqref{we1} based on Lemma \ref{lemma31} and the fixed-point technique, following \cite{SR} or \cite{LZ23}. We omit the proof here. 

\begin{thm}[Local well-posedness]\label{localwp} Assume that $b_1, b_2$ satisfy condition \eqref{we1b}. Given a positive constant $\delta,$ there exists a unique solution $(u,v)\in C^{2+\alpha,1+\alpha/2}(\bar{Q}) \times C^{2+\alpha,1+\alpha/2}(\bar{Q})$ for any $(\tilde{f},\tilde{g}) \in V_{\delta}\times \bar{V}_{\delta},$ where
\[V_{\delta}=\{\tilde{f}\in C^{2+\alpha}(\bar{\Omega}):  \tilde{f}= \tilde{f}_{x_i}=\tilde{f}_{x_i x_j}=0,\, i,j=1,\cdots, n\  \text{on} \  \Gamma ,\text{and}\  \Vert \tilde{f} \Vert_{C^{2+\alpha}(\bar{\Omega})} \leq \delta\}, \]
\[\bar{V}_{\delta}=\{\tilde{g}\in C^{2+\alpha}(\bar{\Omega}):  \tilde{g}= \tilde{g}_{x_i}=\tilde{g}_{x_i x_j}=0,\, i,j=1,\cdots, n\  \text{on} \  \Gamma ,\text{and}\  \Vert \tilde{g} \Vert_{C^{2+\alpha}(\bar{\Omega})} \leq \delta\}.\]
\end{thm}

Next, regarding nonlinearity, we state the global well-posedness of strong solutions to the semilinear parabolic equation below:
\begin{equation}\label{wea}
\begin{cases} 
\partial_{t}u-d_1\Delta u+a_1(x,t,u,v)= 0, &\  \text{in}\   Q,\\ 
\partial_{t}v-d_2\Delta v+a_2(x,t,u,v)= 0, &\  \text{in}\   Q,\\ 
\partial_{\nu}u=\partial_{\nu}v=0, &\  \text{on} \  \Sigma,\\ 
u(x,0)=f(x), &\  \text{in}\   \Omega,\\
v(x,0)=g(x), &\  \text{in}\   \Omega,\\  
\end{cases} 
\end{equation} 
where $f,g \in H^{1}_{0}(\Omega),$ $a_j: Q\times \mathbb{R} \rightarrow \mathbb{R}$ $(j=1,2)$ satisfy the growth conditions:
\begin{equation}\label{wea1}
\mathop{\text{lim} \, \text{sup}}_{u\rightarrow \infty}\frac{\partial_{u}a_j(x,t,u,v)}{\text{ln}^{1/2}\vert u \vert}\  \text{and} \   \mathop{\text{lim} \, \text{sup}}_{v\rightarrow \infty}\frac{\partial_{v}a_j(x,t,u,v)}{\text{ln}^{1/2}\vert v \vert}, \  \text{uniformly for} \  (x,t)\in Q.\end{equation}  
Moreover, $a_j$ also satisfy
\begin{equation}\label{wea2}
a_j(\cdot, \cdot, 0,0)\in L^2(Q), \  a_j(\cdot, \cdot, u_0,v_0)\in L^2(Q),\  a_j(x,t,\cdot, \cdot)\in C^1(\mathbb{R}).\end{equation}

Then, the global well-posedness result of \eqref{wea} can be given as:
\begin{thm}[Global well-posedness\cite{LL23}]\label{globalwp}
 Assume that $a_j$ satisfies conditions \eqref{wea1} and \eqref{wea2}. Then for any $f$ and $g$ belonging to $ H^{1}_{0}(\Omega),$ the semilinear parabolic equation \eqref{wea} admits a unique strong solution $(u,v) \in H^{2,1}(Q)\times H^{2,1}(Q) .$
 \end{thm}

{\centering \section{ANALYSIS OF HIGH-ORDER VARIATION}  \label{highva}}
In this part, we introduce the high-order variation method coupled with successive linearization, which is a novel approach for us to ensure the positivity of the solutions with physical background and apply it to the diffusive Bazykin model. 
We can see from the comparison below that this approach improves on the classical method of high-order variation, by bringing us the same result in the first-order linearization and a simpler form for the second-order linearization.

\subsection{High-order linearization at $(u_0,v_0)$}  \label{highlinear}
To better understand the advantages of the high-order variation method, we first review the application of high-order linearization method into the model \eqref{modac}.   Let $F(x,t,u,v) \in \mathcal{A},$ $G(x,t,u,v) \in \mathcal{B}.$   Let $ f(x; \varepsilon )=u_0+\sum\limits^{N}_{i=1} \varepsilon_{i}f_{i},$ $ g(x; \varepsilon )=v_0+\sum\limits^{N}_{i=1} \varepsilon_{i}g_{i}$ in $\Omega,$ where $f_{i}, g_{i} \in C^{2+\alpha}(\mathbb{R}^{n})$ and $\varepsilon=(\varepsilon_{1},\dots, \varepsilon_{N}) \in \mathbb{R}^{N} $ as $\vert \varepsilon \vert=\vert \varepsilon_{1}\vert + \cdots +\vert \epsilon_{N}\vert $ small enough. One characteristic of this high-order linearization method worth noticing is that when $u_0$ equals 0, $f(x;\varepsilon)$ is not non-negative for sure since the input values $f_i$ may vary around $0,$ which can happen in mathematical analysis but induce problems for physical studies. The same problem arises for $v_0.$ Each $f_i$ and $g_i$ matters in this linearized method.

From the previous section, we know that there exists a unique solution $(u(x,t;\varepsilon), v(x,t;\varepsilon))$ of \eqref{modac}. If $\varepsilon=0,$ then by the admissibility conditions on $F$ and $G,$ we have $( \tilde{u} , \tilde{v}):=(u(x,t;0), v(x,t;0)) =(u_0,v_0)$ as a fixed initial value pair. Let 
\begin{equation}u^{(1)}:=\partial_{\varepsilon_1}u\vert_{\varepsilon=0}\  \text{where}\  \mathop{\text{lim}}\limits_{\vert \varepsilon\vert \to 0}\frac{\vert u(x,t;\varepsilon)- u(x,t;0)-u^{(1)}(x,t)\varepsilon \vert}{\vert \varepsilon_1 \vert}=0 \quad \text{in}\ Q,\end{equation}
\begin{equation}v^{(1)}:=\partial_{\varepsilon_1}v\vert_{\varepsilon=0}\  \text{where}\  \mathop{\text{lim}}\limits_{\vert\varepsilon \vert\to 0}\frac{ \vert v(x,t;\varepsilon)- v(x,t;0)-v^{(1)}(x,t)\varepsilon\vert}{\vert \varepsilon_1 \vert}=0\quad \text{in}\ Q.\end{equation}

Then, in the method of high-order linearization,  a new system is built with the unknowns $( u^{(1)}, v^{(1)}).$ This is possible by Theorems \ref{localwp} and \ref{globalwp}, since the assumptions are satisfied for $F\in \mathcal{A}$ and $G\in \mathcal{B}.$
It is also known from the definitions of the admissible sets $\mathcal{A}$ and $\mathcal{B}$  that both $F^{(1)}_{u}(u,v)=F^{(1)}_{v}(u,v)=0$ and $G^{(1)}_{u}(u,v)=0$ hold true at all times. Therefore, we have that $( u^{(1)}, v^{(1)})$  satisfies the following system:
\begin{equation}\label{ifl}
\begin{cases} 
\partial_t u^{(1)}(x,t)-d_1\Delta u^{(1)}(x,t)= 0, &\  \text{in}\   Q,\\ 
\partial_t v^{(1)}(x,t)-d_2\Delta v^{(1)}(x,t)= G^{(1)}_{v}v^{(1)}(x,t), &\  \text{in}\   Q,\\ 
\partial_{\nu} u^{(1)}(x,t)=\partial_{\nu} v^{(1)}(x,t)=0, &\  \text{on}\   \Sigma,\\ 
u^{(1)} (x,0)= f_1(x) , \; v^{(1)} (x,0)= g_1(x)&\  \text{in}\  \Omega. 
\end{cases} 
\end{equation} 
 
A similar first-order system linearized around other $\epsilon_l$'s can be constructed by defining the following:
\[u^{(l)}:=\partial_{ \varepsilon_{l}}u\vert_{\varepsilon=0}, \  l\in \{1,2,\dots,N\},\] 
\[v^{(l)}:=\partial_{\varepsilon_{l}}v\vert_{\varepsilon=0}, \  l\in \{1,2,\dots,N\}.\] 
\begin{rmk}\label{rm41}
   One thing worth noticing in this method application is that classically we do not restrict $f_1(x)$ and $g_1(x)$ no matter what the values of $u_0$ and $ v_0$ are. Hence, we cannot ensure the positivity of $u^{(1)}(x,t)$ and $v^{(1)}(x,t).$  
\end{rmk}

Next, we consider the second-order linearization, which we define as 
\[u^{(1,2)}:=\partial_{\varepsilon_{1}}\partial_{\varepsilon_{2}} u\vert_{\varepsilon=0}, \  v^{(1,2)}:=\partial_{\varepsilon_{1}}\partial_{\varepsilon_{2}} v\vert_{\varepsilon=0}.\]

The second-order linearization is:
\begin{equation} \label{isl}
\begin{cases} 
u^{(1,2)}_t-d_1\Delta u^{(1,2)} = F^{(1,2)}_{vu}v^{(1)}u^{(2)}+F^{(1,2)}_{uv}u^{(1)}v^{(2)} +&\  \\
\qquad  \qquad  \qquad  \qquad F^{(1,2)}_{vv}v^{(1)}v^{(2)}+F^{(1,2)}_{uu}u^{(1)}u^{(2)} , &\  \text{in} \   Q,\\ 
v^{(1,2)}_t-d_2\Delta v^{(1,2)}= G^{(1,2)}_{vu}v^{(1)}u^{(2)}+G^{(1,2)}_{uv}u^{(1)}v^{(2)} +&\  \\
\qquad  \qquad  \qquad  \qquad G^{(1,2)}_{vv}v^{(1)}v^{(2)}+G^{(1,2)}_{uu}u^{(1)}u^{(2)}+  G^{(1)}_{v} v^{(1,2)} , &\  \text{in} \   Q,\\ 
\partial_{\nu} u^{(1,2)}(x,t)=\partial_{\nu} v^{(1,2)}(x,t)=0, &\  \text{on}\   \Sigma,\\ 
u^{(1,2)} (x,0)=0, \; v^{(1,2)} (x,0)= 0,&\  \text{in}\  \Omega.
\end{cases} 
\end{equation} 

It is apparent that the second-order linearized system \eqref{isl} depends on the first-order linearized system \eqref{ifl}. What we achieved after studying \eqref{ifl} would become known in \eqref{isl}, reflecting this method's successive property.

In general, since $ f(x)$ and $ g(x)$ are infinitely differentiable, if we use $\{f_1,f_2,\dots,f_{N}\} $ and $\{g_1,g_2,\dots,g_{N}\}$  in $C^{2+\alpha}$ to expand \eqref{modac}, by Theorems \ref{localwp} and \ref{globalwp}, $u$ and $v$ are infinitely differentiable, so
we can get an $N$-th order linearized system by defining:
\[u^{(l_1,l_2,\dots,l_{N})}(x,t)=\partial_{\varepsilon_{l_{N}}} \cdots\partial_{\varepsilon_{l_2}}\partial_{\varepsilon_{l_1}} u\vert_{\varepsilon=0},\ l_1,l_2,\dots,l_N\in \{1,2,\dots,N\},\] 
\[v^{(l_1,l_2,\dots,l_{N})}(x,t)=\partial_{\varepsilon_{l_{N}}} \cdots\partial_{\varepsilon_{l_2}}\partial_{\varepsilon_{l_1}} v\vert_{\varepsilon=0},\ l_1,l_2,\dots,l_N\in \{1,2,\dots,N\}.\] 

The core of this method is that the solutions of lower-order terms decide the nonlinear terms in higher-order systems. We can use mathematical induction to prove Theorem \ref{mainthm} as \cite{SR,LZ22,LMZ23,LZ23} did. However, in this paper, we apply the following method in the main proof.

\subsection{High-order variation at $(u_0,v_0)$} \label{highvary}
As we discussed above, the classical method of high-order linearization may render the results physically meaningless, since the solutions to the model problem \eqref{modac} are not non-negative. To  address this non-negativity constraint on the species population densities $u$ and $v,$ we introduce the comparatively new technique, the high-order variation method. The biggest difference from the former method is how we define the initial input, for a positive-chosen initial value, to ensure the non-negativity for the solutions. Furthermore, it turns out that this method can also offer a simpler form for proving the uniqueness in higher orders.

Consider the system \eqref{modac}. If $(u_0,v_0)$ is known as a solution to \eqref{modac}, let
\[f(x;\varepsilon)=u_0+\sum\limits^{N}_{i=1} \varepsilon^i f_{i},\  g(x;\varepsilon)=v_0+\sum\limits^{N}_{i=1} \varepsilon^i g_{i} \quad \text{in}\ \Omega,\]
where $f_{i}, g_{i} \in C^{2+\alpha}(\Omega),$  $ \varepsilon \in \mathbb{R}^{+}.$ When $u_0=0,$ we ask $f_{1}\geq 0,$ and when $v_0=0,$ we ask $g_{1} \geq 0.$ Hence $f,g \geq 0$ in $\Omega$ as $\varepsilon$ goes to $0.$ There is no need to add restrictions for $f_{1}$  if $u_0$ is known to be positive constant, this holds similarly for $g_{1}$ and its corresponding $v_0.$ This ensures the positivity of $u,v.$

Now we define the first-order variation form as:
\[u^{(I)}:=\partial_{ \varepsilon}u\vert_{\varepsilon=0}\  \text{where} \  \mathop{\text{lim}}\limits_{\vert\varepsilon \vert \to 0}\frac{\vert u(x,t;\epsilon)- u(x,t;0)-u^{ (I) }(x,t)\varepsilon \vert}{\vert \varepsilon \vert}=0\quad  \text{in} \ Q,\] 
\[v^{(I)}:=\partial_{\varepsilon}v\vert_{\varepsilon=0}\  \text{where} \  \mathop{\text{lim}}\limits_{\vert\varepsilon \vert \to 0}\frac{\vert v(x,t;\epsilon)- v(x,t;0)-v^{ (I) }(x,t)\varepsilon \vert}{\vert \varepsilon \vert}=0\quad \text{in} \ Q. \] 

The first-order linearization would satisfy: 
\begin{equation} \label{ifv}
\begin{cases} 
u^{ (I) }_{t}(x,t)-d_1\Delta u^{(I) }(x,t)=0, &\  \text{in}\   Q,\\ 
v^{(\uppercase\expandafter{\romannumeral1}) }_{t}(x,t)-d_2\Delta v^{(I) }(x,t)=G^{(I) }_{v} v^{(I)}, &\  \text{in}\    Q,\\ 
\partial_{\nu} u^{(I)}(x,t)=\partial_{\nu} v^{(I)}(x,t)=0, &\  \text{on}\   \Sigma,\\
u^{(I)} (x,0)=f_1(x), \; v^{(I)} (x,0)=  g_1(x),& \  \text{in}\    \Omega.
\end{cases} 
\end{equation} 

Since the definition of the first-order variation form is the same as that for high-order linearization, \eqref{ifv} keeps the same form with \eqref{ifl}. However, the non-negative input of the first-order variation system ensures its non-negativity, which makes our results more physically realistic.

The second-order variation is defined as
 \[ u^{(II)}:=\partial^{2}_{\varepsilon}u\vert_{\varepsilon=0}, \, v^{(II)}:=\partial^{2}_{\varepsilon}v\vert_{\varepsilon=0}. \]
 
Then the second-order variation system is given as 
\begin{equation} \label{isv}
\begin{cases} 
u^{(II)}_t-d_1\Delta u^{(II)}  = (F^{(II)}_{vu}+F^{(II)}_{uv})v^{(I)}u^{(I)}+&\ \\
\qquad  \qquad  \qquad  \qquad F^{(II)}_{uu}(u^{(I)})^2+F^{(II)}_{vv}(v^{(I)})^2 , &\  \text{in}\    Q,\\ 
v^{(II)}_t-d_2\Delta v^{(II)} -G^{(I)}_{v} v^{(II)}=( G^{(II)}_{vu}+G^{(II)}_{uv})v^{(I)}u^{(I)}+&\ \\
\qquad  \qquad  \qquad  \qquad G^{(II)}_{vv}(v^{(I)})^2 +G^{(II)}_{uu} (u^{(I)})^2, &\  \text{in}\    Q,\\ 
\partial_{\nu} u^{(II)}(x,t)=\partial_{\nu} v^{(II)}(x,t)=0, &\  \text{on}\   \Sigma,\\
u^{(II)} (x,0)=2f_2, \; v^{(II)} (x,0)= 2g_2,& \  \text{in} \   \Omega. 
\end{cases} 
\end{equation} 

In this system, $f_2$ and $g_2$ are given arbitrarily no matter what the initial value $u_0, v_0$ is, since the positivity of $u,v$ is guaranteed by the positivity of $f_1,g_1.$

Similarly, for $k\in {2, \cdots, N},$ we consider 
\[u^{( k)}=\partial^{k}_{\varepsilon}u|_{\varepsilon=0},\  \  v^{( k)}=\partial^{k}_{\varepsilon}v|_{\varepsilon=0},\]
and we can generate a series of parabolic systems to determine the high-order Taylor coefficients of $F$ and $G.$ Via this definition, the system at each order has a simpler form than those obtained by the high-order linearization method. Moreover, the solutions of lower-order terms successively decide the non-linear terms in higher-order systems, which is a common feature of both methods.

\subsection{Advantages of choosing high-order variation method}
Comparing sections \ref{highlinear} and \ref{highvary}, we can find two obvious advantages of high-order variation method over high-order linearization method.

First and most importantly, the high-order variation method has better behavior when we study the inverse problem of a physical model  requiring positive solutions (or larger than a fixed number). We ask $f_1, g_1>0$ for $u_0=0,v_0=0,$ which gives $u^{(I)}(x,t), v^{(I)}(x,t)>0$ from the maximum principle of parabolic equations applied to \eqref{ifv}, when the first order derivatives of $G$ fulfil coercivity conditions. It gives meaning for physical uses and is also highly applicable to a wide range of PDEs.

Moreover, we define a simpler second-order variation form for the system, which helps significantly reduce the calculation in later proofs. To be specific, to solve \eqref{isv}, we only need to solve \eqref{ifv} as a premise. On the other hand, if we want to solve \eqref{isl}, we are required to solve two first-order linearization systems, namely \eqref{ifl} and a similar system for  $(u^{(2)},v^{(2)}).$ High-order variation method helps us save the calculation process.

Here, we also stress that we made a massive improvement to the method of high-order variation compared with existing articles.\cite{LZ23} only used this new technique around the solution $(0,1),$ \cite{LL23, DL23, DL23mfg} applied high-order variation around a pair of trivial solutions $(0,0),$ and in this paper, we allow for more general solutions $(u_0,v_0)$ of the system, yet maintain the positivity of $u$ and $v$ in a general form. Hence, our results cover the cases of applying high-order variation around $(1,0)$ and $(0,0)$ previously used, and this technical development shall help tackle a broader range of nonlinear inverse problems related to physical systems especially those with positivity constraint.

{\centering \section{PROOF OF THEOREM 2.4}   \label{proofm}}

We mainly apply the method of high-order variation to prove Theorem \ref{mainthm} in this section, and we begin with two crucial auxiliary lemmas.

\begin{lem}\label{formu}
Consider 
\begin{equation}\label{lem1}
\begin{cases} 
\partial_t u(x,t)-p\Delta u(x,t)= 0, &\  \text{in} \   Q,\\ 
\partial_{\nu} u(x,t)=0, &\  \text{on} \   \Sigma,
\end{cases} 
\end{equation}
where $p$ is a constant. There exists a sequence of solutions $u(x,t)$ to the system \eqref{lem1} such that $u(x,t)=e^{\lambda t} h(x;\lambda)$ for some $\lambda \in \mathbb{R}^n$ and $h(x;\lambda) \in C^{2}(\Omega).$ In particular, $h(x;\lambda)$ is not necessarily 0, and $\frac{\lambda}{p}$ is its corresponding eigenvalue.\end{lem}

\textit{Proof.} Let $\frac{\lambda}{p}$ be an Neumann Laplacian eigenvalue and $h(x;\lambda)$ be its eigenfunction.
\begin{equation}\label{lem1p}
\begin{cases} 
-\Delta h(x;\lambda)= \frac{\lambda}{p}h(x;\lambda), &\  \text{in}\    Q,\\ 
\partial_{\nu} h(x;\lambda)=0, &\  \text{on} \   \Sigma.
\end{cases} 
\notag
\end{equation}
It is obvious to see that $ u(x,t)=e^{\lambda t} h(x;\lambda) $ is a solution of \eqref{lem1}.  $ \hfill{\square}$
\begin{lem}\label{formv}
Consider 
\begin{equation}\label{lem2}
\begin{cases} 
\partial_t v(x,t)-q\Delta v(x,t)+k v(x,t)= 0, &\  \text{in} \   Q,\\ 
\partial_{\nu} v(x,t)=0, &\  \text{on} \   \Sigma,
\end{cases} 
\end{equation}
where $q$ and $k$ are constants. There exists a sequence of solutions $v(x,t)$ to the system \eqref{lem2} such that $v(x,t)=e^{\mu t} l(x;\mu)$ for some $\mu \in \mathbb{R}^n$ and $l(x;\lambda) \in C^{2}(\Omega).$ In particular, $l(x;\mu)$ is not necessarily 0, and $\frac{\mu}{q}$ is its corresponding eigenvalue.\end{lem}

\textit{Proof.}  Since $k$ is a constant, by letting $v(x,t)=u(x,t)e^{-kt},$ we can transform \eqref{lem2} into
\begin{equation}\label{lem2p}
\begin{cases} 
\partial_t u(x,t)-q\Delta u(x,t)= 0, &\  \text{in} \   Q,\\ 
\partial_{\nu} u(x,t)=0, &\  \text{on} \   \Sigma,
\end{cases} 
\notag
\end{equation}
which is the case in Lemma \ref{formu}. $ \hfill{\square}$

Now we present the proof of Theorem \ref{mainthm} with the high-order variation method. Let 
\[f(x)=u_0+\sum\limits^{N}_{i=1}\varepsilon^i f_i,\ g(x)=v_0+\sum\limits^{N}_{i=1}\varepsilon^i g_i.\]
The innovative crux of the high-order variation method requires asking $f_1\geq 0$ if $u_0=0$ as the initial value, and similarly $g_1\geq 0$ if $v_0=0.$ Through this setting, we can always ensure the non-negativity for $u$ and $v.$ 

\subsection{Unique recovery of the first-order coefficient}

For $j=1,2,$ we consider 
\begin{equation} \label{genc}
\begin{cases} 
\partial_t u_{j}(x,t)-d_1\Delta u_j(x,t)=   F_j(x,t,u,v) &\  \text{in} \    Q,\\ 
\partial_t v_{j}(x,t)-d_2\Delta v_j (x,t)= G_j(x,t,u,v)&\  \text{in} \    Q,\\ 
\partial_{\nu}u_{j}(x,t)=\partial_{\nu}v_{j}(x,t)=0 &\  \text{on}\    \Sigma,\\ 
u_j (x,0)=f(x)  &\  \text{in} \   \Omega,\\ 
v_j (x,0)=g(x) & \  \text{in}\    \Omega. 
\end{cases} 
\end{equation}
The system \eqref{genc} is known to have a solution $(u_0,v_0),$ when $F_j \in \mathcal{A},$ and $G_j \in \mathcal{B}.$ 

The first step of the unique identifiability problem is to prove the equality of the first-order Taylor coefficients  of $F_j$ and $G_j$ for $j=1,2.$ Since $F_j \in \mathcal{A} (j=1,2),$ we know by definition that:
\begin{equation} F^{(1)}_{u,1}=F^{(1)}_{u,2}=F^{(1)}_{v,1}=F^{(1)}_{v,2}=G^{(1)}_{u,1}=G^{(1)}_{u,2}=0.\notag \end{equation}
Therefore, it only remains to show $G^{(1)}_{v,1}=G^{(1)}_{v,2}.$ Meanwhile, we define the first-order variation for $u$ and $v$ as introduced in section \ref{highvary}.

Take $f_1=1.$ Then we have 
\begin{equation} \label{1stu}
\begin{cases} 
\partial_{t} u^{ (I) }_{j}(x,t)-d_1\Delta u^{(I) }_{j}(x,t)=0 &\  \text{in}\   Q,\\ 
\partial_{t} v^{ (I) }_{j}(x,t)-d_2\Delta v^{(I) }_{j}(x,t)=G^{(I) }_{v,j} v^{(I)}_{j}&\  \text{in}\    Q,\\ 
\partial_{\nu} u^{(I)}_{j}(x,t)=\partial_{\nu} v^{(I)}_{j}(x,t)=0 &\  \text{on}\   \Sigma,\\
u^{(I)}_{j} (x,0)=1 & \  \text{in}\    \Omega, \\ 
v^{(I)} _{j}(x,0)=  g_1(x)\geq  0 & \  \text{in}\    \Omega.
\end{cases} 
\end{equation} 

Let  $\bar{u}(x,t):=u^{(I)}_1 (x,t)-u^{(I)}_2(x,t).$ When $\mathcal{M}^{+}_{F_1,G_1}=\mathcal{M}^{+}_{F_2,G_2}$ is known as a precondition, $\bar{u}(x,t)$ satisfies
\begin{equation} \label{baru1}
\begin{cases} 
\partial_t \bar{u}(x,t)-d_1\Delta  \bar{u}(x,t) = 0, &\  \text{in} \   Q,\\ 
\partial_{\nu}  \bar{u}(x,t)= \bar{u}(x,t)=0, &\  \text{on} \   \Sigma,\\ 
 \bar{u}(x,0)=\bar{u}(x,T)=0,& \  \text{in}\    \Omega. 
\end{cases} 
\end{equation} 
Since this is simply the heat equation, it is easy to see $\bar{u}=0$ and thus $u^{(I)}_1 (x,t)=u^{(I)}_2(x,t):=u^{(I)}(x,t).$ Meanwhile, by Lemma \ref{formu}, there exists $\lambda \in \mathbb{R}$ and $h(x)\in C^{2}(\Omega)$ such that
\begin{equation} \label{splitu}
u^{(I)}(x,t)=e^{\lambda t}h(x) \end{equation}
satisfies equations for $u^{(I)}_{j}\, (j=1,2)$ in  \eqref{1stu}. 

On the other hand, let  $\bar{v}(x,t)=v^{(I)}_1 (x,t)-v^{(I)}_2(x,t).$ Since $\mathcal{M}^{+}_{F_1,G_1}=\mathcal{M}^{+}_{F_2,G_2}$ is known as a precondition, $\bar{v}(x,t)$ satisfies
\begin{equation}\label{1stv}
\begin{cases} 
\partial_t \bar{v}-d_2\Delta  \bar{v}=  G^{(I) }_{v,1} \bar{v}+(G^{(I) }_{v,1}-G^{(I) }_{v,2})v^{(I)}_2 &\  \text{in} \   Q,\\ 
\partial_{\nu}  \bar{v}(x,t)= \bar{v}(x,t)=0 &\  \text{on} \   \Sigma,\\ 
 \bar{v}(x,0)=\bar{v}(x,T)=0& \  \text{in}\    \Omega. 
\end{cases} 
\end{equation} 

Let $\omega$ be a solution of the following system 
\begin{equation}\label{omegaform1}
-\partial_{t }\omega -d_2 \Delta \omega -G^{(I) }_{v,1} \omega=0 \  \text{in} \  Q,
\end{equation}
where $G^{(I) }_{v,1}$ is an unknown constant, and the CGO solution to $\omega$ is easy to seek from \eqref{omegaform1} as
\begin{equation}\label{CGOfw}
\omega=e^{(|\xi|^2-G^{(I) }_{v,1})t-\frac{i}{\sqrt{d_2}}\xi \cdot x},
\end{equation}
with $i=\sqrt{-1}$ for $\xi \in \mathbb{R}^n.$

Then we multiply $\omega$ on both sides of \eqref{1stv} and conduct integration by parts to achieve
\begin{equation}\label{simGv}\int_{Q} (G^{(I) }_{v,1}-G^{(I) }_{v,2})v^{(I)}_2 \omega dxdt=0. \end{equation}

Since $v^{(I)}_{j}$ is of the form \eqref{lem2} with $d_2, G^{(I)}_{v,j}$ constant, by Lemma \ref{formv}, there exist $\mu \in \mathbb{R}$ and $l(x)\in C^{\infty}(\Omega)$ such that $e^{\mu t}l(x)$ satisfies equations for $v^{(I)}_{j}\, (j=1,2)$ in \eqref{1stu}. By the uniqueness of the solution of second-order parabolic equations, we have
\begin{equation}\label{nearfv} v^{(I)}_{2}(x,t)=e^{\mu t}l(x).\end{equation}

Substituting \eqref{CGOfw} and \eqref{nearfv} into \eqref{simGv}, we know \eqref{simGv} satisfies:
\begin{equation}\label{subsdirec1}
    \int^{T}_{0}e^{\mu t}e^{(|\xi|^2-G^{(I) }_{v,1})t}dt \int_{\Omega} (G^{(I) }_{v,1}-G^{(I) }_{v,2})l(x;\mu) e^{-\frac{i}{\sqrt{d_2}}\xi \cdot x} dx=0, 
\notag\end{equation}
which yields 
\begin{equation}\label{variableapart1}
\int_{\Omega} (G^{(I) }_{v,1}-G^{(I) }_{v,2})l(x;\mu) e^{-\frac{i}{\sqrt{d_2}}\xi \cdot x} dx=0.
\notag\end{equation}
Since this holds for any Neumann eigenfunction $l(x;\mu)$ of $\Delta,$ we obtain
\[G^{(I) }_{v,1}=G^{(I) }_{v,2}.\]
 
\begin{rmk}\label{solv}
\textit{One thing worth mentioning here is that we choose} $f_1(x,t)=1, g_1(x,t)\geq 0.$ \textit{ The positivity of } $f_1$ \textit{means that for all} $\varepsilon, u(x;\varepsilon)>0$ \textit{in} $\Omega.$ \textit{By the maximum principle for the heat problem,} $u(x,t)>0$ \textit{in} $Q.$ Then we can solve the equation \eqref{1stu} for $v^{(I)}_j(x,t)$ and obtain the unique solution: 
\begin{equation}\label{exprv}v^{(I)}_j(x,t)=\int_{\Omega}\Psi(x-y,t)g_1(y)dy,\end{equation}
\textit{where} $\Psi(x,t)$ \textit{is the Green's function operator for} $\partial_t-d_2\Delta -G^{(I) }_{v}$ \textit{on the bounded domain} $\Omega.$ \textit{To ensure} $v^{(I)}_j(x,t) \geq 0,$ \textit{we can ask} $G^{(I) }_{v} \leq 0.$ \textit{High-order variation ensures the physical meaning of} $u,v.$
\end{rmk}

\begin{rmk}\label{1stnonneg}
\textit{In this proof, we choose a concrete} $f_1(x,t)=1$ \textit{and a comparatively less specific initial function } $g_1(x,t)\geq 0.$ \textit{We can also relax the requirement of} $f_1$ \textit{into any non-negative function. All the assumptions on initial functions serve for the convenience of the proof.}
\end{rmk}

\subsection{Unique recovery of the second-order coefficients }

In this section, we shall see the advantage of using high-order variation method to simplify the calculation process. However, the recovery for coefficients in second-order variation systems is much more complicated than those in first-order variation systems.  We divide all unknown coefficients into two sets, $A:=\{F^{(II)}_{uu},F^{(II)}_{uv},F^{(II)}_{vu},F^{(II)}_{vv}\}$ and $B:=\{G^{(II)}_{uu},G^{(II)}_{uv},G^{(II)}_{vu},G^{(II)}_{vv}\}.$
We will recover one coefficient for each set, assuming the other three are known. Then $A, B$ generate 16 combinations for us to realize the recovery mission.

We start by recalling the second-order variation system as described in the previous section, given by \eqref{isv}. Let
\begin{equation} u^{(II)}_j:=\partial^{2}_{\varepsilon}u_j\vert_{\varepsilon=0}, \  j=1,2.  \notag \end{equation}
\begin{equation}  v^{(II)}_j:=\partial^{2}_{\varepsilon}v_j\vert_{\varepsilon=0}, \  j=1,2. 
\notag \end{equation}

The second-order variation system follows:
\begin{equation}  \label{secgs}
\begin{cases} 
\partial_t u^{(II)}_j-d_1\Delta u^{(II)}_j  = (F^{(II)}_{vu}+F^{(II)}_{uv})v^{(I)}u^{(I)}+&\ \\
\qquad  \qquad  \qquad  \qquad F^{(II)}_{uu,j} (u^{(I)})^2+F^{(II)}_{vv}(v^{(I)})^2 , &\  \text{in}\    Q,\\ 
\partial_t  v^{(II)}_j-d_2\Delta v^{(II)}_j -G^{(I)}_{v} v^{(II)}_j=( G^{(II)}_{vu}+G^{(II)}_{uv})v^{(I)}u^{(I)}+&\ \\
\qquad  \qquad  \qquad  \qquad G^{(II)}_{vv}(v^{(I)})^2 +G^{(II)}_{uu,j} (u^{(I)})^2  , &\  \text{in}\    Q,\\ 
\partial_{\nu} u^{(II)}_j(x,t)=\partial_{\nu} v^{(II)}_j(x,t)=0, &\  \text{on}\   \Sigma,\\
u^{(II)} (x,0)=2f_2, \; v^{(II)} (x,0)= 2g_2,& \  \text{in} \   \Omega. 
\end{cases} 
\end{equation} 
Here we have united the notation of $u^{(I)}_j, v^{(I)}_j,$ $j=1,2$ into $u^{(I)}$ and $v^{(I)}$ respectively. Indeed, this is possible as follows: since it is obvious from the uniqueness of solution of the heat equation \eqref{baru1} that $\bar{u}=0,$ therefore $u^{(I)}_{1}=u^{(I)}_{2}.$ On the other hand, we have already recovered $G^{(I) }_{v,1}=G^{(I) }_{v,2}$ in the previous subsection, and from  expression \eqref{exprv}, thus $v^{(I)}_{1}=v^{(I)}_{2}.$

The first-order variation term $u^{(I)}, v^{(I)}$ in \eqref{secgs} depends on the solution of \eqref{1stu} for arbitrary $f_1, g_1.$ However, instead of the non-negative restrictions for initial data in \eqref{1stu}, $f_2,g_2$ can be chosen arbitrarily, because both $u^{(II)}$ and $v^{(II)}$ do not need to be strictly positive, since they are higher order perturbations of the non-negative lower first-order terms $u^{(I)},$ $v^{(I)}.$ Moreover, compared with using high-order linearization method, we reduce calculations in solving the second-order system. Specifically, we only need $u^{(I)},v^{(I)}$ to solve the second-order variation system  \eqref{secgs}, but we need two first-order linearization systems for $u^{(1)},v^{(1)}$ and $u^{(2)},v^{(2)}$ to solve one second-order linearization system involving $u^{(1,2)},v^{(1,2)}.$

Next, let $\hat{u}(x,t)=u^{(II)}_1 (x,t)-u^{(II)}_2 (x,t).$ Using $\mathcal{M}^{+}_{F_1,G_1}=\mathcal{M}^{+}_{F_2,G_2}$ as a precondition, we have the following equations from \eqref{secgs}:
\begin{equation} \label{hatu2}
\begin{cases} 
\partial_t \hat{u}-d_1\Delta \hat{u}  = (F^{(II)}_{vu,1}-F^{(II)}_{vu,2})v^{(I)}u^{(I)}+(F^{(II)}_{uv,1}-F^{(II)}_{uv,2})v^{(I)}u^{(I)}+&\ \\
\qquad  \qquad  \qquad  \qquad (F^{(II)}_{uu,1}-F^{(II)}_{uu,2} )(u^{(I)})^2+(F^{(II)}_{vv,1}-F^{(II)}_{vv,2})(v^{(I)})^2 , &\  \text{in}\    Q,\\ 
\partial_{\nu} \hat{u}(x,t)=\hat{u}(x,t)=0 , &\  \text{on}\    \Sigma,\\ 
\hat{u}(x,0)=\hat{u}(x,T)=0,  &\  \text{in} \   \Omega.
\end{cases} 
\end{equation}

Let $w$ be the solution of $-\partial_t w-d_1\Delta w=0$ in $Q,$ multiply it on both sides of \eqref{hatu2} and then integrate by parts, we obtain
\begin{equation} \label{hatFu} 
\begin{split}
   \int_{Q} \big[(F^{(II)}_{vu,1}-F^{(II)}_{vu,2})v^{(I)}u^{(I)}+(F^{(II)}_{uv,1}-F^{(II)}_{uv,2})v^{(I)}u^{(I)}+(F^{(II)}_{uu,1}-F^{(II)}_{uu,2}) (u^{(I)})^2\\
   +(F^{(II)}_{vv,1}-F^{(II)}_{vv,2}) (v^{(I)})^2 \big] w dxdt=0. 
\end{split}
\end{equation}

Combining with \eqref{splitu}, \eqref{nearfv} and substituting $w=e^{|\xi|^2 t- \frac{i}{\sqrt{d_1}}\xi \cdot x}$ into \eqref{hatFu}, we obtain four cases of recovering $F^{(II)}$ given different assumptions.

\textbf{\textit{Recovery of $F^{(II) }_{vu,j}$ and $F^{(II) }_{uv,j}$} }
The recovery of  $F^{(II)}_{vu,j}$ and $F^{(II)}_{uv,j}$ are very similar, so we only illustrate the unique identifiability procedure for $F^{(II) }_{vu,j}.$ 

We assume  $F^{(II) }_{uu},F^{(II) }_{uv},F^{(II) }_{vv}$  are known, then \eqref{hatFu} becomes
\begin{equation} \label{Fvu2}
\int_{\Omega} (F^{(II)}_{vu,1}-F^{(II)}_{vu,2})h(x;\lambda)l(x;\mu) e^{- \frac{i}{\sqrt{d_1}}\xi \cdot x} dx=0. 
\end{equation}

Picking $h(x;\lambda),$ $l(x;\mu)$ to be any non-zero eigenfunction of $\Delta,$ we have from \eqref{Fvu2} that
 \[F^{(II) }_{vu,1}=F^{(II) }_{vu,2}.\]
 
 We can obtain the uniqueness of $F^{(II) }_{uv}$ similarly.
 
\textbf{\textit{Recovery of $F^{(II) }_{uu,j}$} } Now we assume $F^{(II) }_{vu},F^{(II) }_{uv},F^{(II) }_{vv}$  are known, and we aim to depict the identity for $F^{(II) }_{uu,j},\ j=1,2.$

From \eqref{hatFu} , we now have
\begin{equation} \label{Fuu2}
\int_{\Omega} (F^{(II)}_{uu,1}-F^{(II)}_{uu,2})h^2(x;\lambda) e^{- \frac{i}{\sqrt{d_1}}\xi \cdot x} dx=0. 
\end{equation}

Since $h(x;\lambda)$ can be chosen to be a non-zero eigenfunction of $\Delta,$ we have
 \[F^{(II) }_{uu,1}=F^{(II)}_{uu,2}.\]

\textbf{\textit{Recovery of $F^{(II) }_{vv,j}$} } Similarly, we can recover $F^{(II) }_{vv,j}$ given $F^{(II) }_{vu},F^{(II) }_{uv}$ and $F^{(II) }_{uu}.$ 

Substituting the given information into \eqref{hatFu} and we have, 
\begin{equation} \label{Fvv2}
\int_{Q} (F^{(II)}_{vv,1}-F^{(II)}_{vv,2})  l^2(x;\mu)  e^{-\frac{i}{\sqrt{d_1}}\xi \cdot x} dx=0. 
\end{equation}

This time, since $l(x;\mu)$ can be chosen to be a non-zero eigenfunction of $\Delta,$ we obtain from \eqref{Fvv2} that
 \[F^{(II) }_{vv,1}=F^{(II) }_{vv,2}.\]

Next, we consider the recovery for coefficients $G^{(II)}.$ One thing worth mentioning here is that we already recovered $G^{(I)}_u$ and $G^{(I)}_v$ in the previous subsection. Let $\hat{v}(x,t)=v^{(II)}_1 (x,t)-v^{(II)}_2 (x,t).$ By $\mathcal{M}^{+}_{F_1,G_1}=\mathcal{M}^{+}_{F_2,G_2},$
we have the following system from \eqref{secgs}:
\begin{equation} \label{hatv2}
\begin{cases} 
\partial_t \hat{v}-d_2\Delta \hat{v} =( G^{(II)}_{vu,1}-G^{(II)}_{vu,2})v^{(I)}u^{(I)}+(G^{(II)}_{uv,1}-G^{(II)}_{uv,2})v^{(I)}u^{(I)}+&\ \\
\qquad  \qquad  \qquad (G^{(II)}_{vv,1}-G^{(II)}_{vv,2})(v^{(I)})^2 +(G^{(II)}_{uu,1}-G^{(II)}_{uu,2})(u^{(I)})^2 +G^{(I)}_{v} \hat{v}, &\  \text{in}\    Q,\\ 
\partial_{\nu} \hat{v}(x,t)=\hat{v}(x,t)=0, &\  \text{on}\   \Sigma,\\
\hat{v}(x,0)=\hat{v}(x,T)= 0,& \  \text{in} \   \Omega. 
\end{cases} 
\end{equation} 

Observe that all the variation terms of $u$ and $v$ involved in \eqref{hatv2} are of a lower variation order, and have already been previously determined. Therefore, our recovery of the coefficients of $G^{(II)}$ is simultaneous.

Let $\omega$ be a solution of the following system 
\begin{equation}\label{omegaforv}  -\partial_{t }\omega -d_2 \Delta \omega -G^{(I) }_{v} \omega=0 \  \text{in} \  Q.
\end{equation} 
Multiplying $\omega$ on both sides of \eqref{hatv2} and integrating by parts, we have:
\begin{equation}\label{hatGv} 
\begin{split}
   \int_{Q} \big[ (G^{(II)}_{vu,1}-G^{(II)}_{vu,2})v^{(I)}u^{(I)}+(G^{(II)}_{uv,1}-G^{(II)}_{uv,2})v^{(I)}u^{(I)}+(G^{(II)}_{uu,1}-G^{(II)}_{uu,2}) (u^{(I)})^2\\
   +(G^{(II)}_{vv,1}-G^{(II)}_{vv,2}) (v^{(I)})^2 \big] \omega dxdt=0. 
\end{split}
\end{equation}
Then, combining \eqref{splitu}, \eqref{nearfv}, substituting the CGO solution $\omega=e^{(|\xi|^2-G^{(I) }_{v,1})t- \frac{i}{\sqrt{d_2}}\xi\cdot x}$ for \eqref{omegaforv}
 into \eqref{hatGv}, we can obtain four cases of recovering $G^{(II)}$ under different assumptions.

\textbf{\textit{Recovery of $G^{(II) }_{vu,j}$ and $G^{(II) }_{uv,j}$} }
The recovery of  $G^{(II)}_{vu,j}$ and $G^{(II)}_{uv,j}$ are very similar, we will only solve the identification problem for $G^{(II) }_{vu,j}$ as an example.

Assume that  $G^{(II) }_{uu},G^{(II) }_{uv},G^{(II) }_{vv}$ are known, then \eqref{hatGv} becomes
\begin{equation} \label{Gvu2}
\int_{\Omega} (G^{(II)}_{vu,1}-G^{(II)}_{vu,2}) l(x;\mu)h(x;\lambda) e^{-\frac{i}{\sqrt{d_2}}\xi\cdot x} dx=0.  
\end{equation}
As in the recovery of  $F^{(II) }_{vu,j}$ and $F^{(II) }_{uv,j},$ for non-zero eigenfunctions $h(x;\lambda),l(x;\mu),$ we obtain $G^{(II)}_{vu,1}=G^{(II)}_{vu,2}.$

 We can obtain the uniqueness of $G^{(II) }_{uv}$ similarly.

\textbf{\textit{Recovery of $G^{(II) }_{uu,j}$} } Now we assume $G^{(II) }_{vu},G^{(II) }_{uv},G^{(II) }_{vv}$ are known, and we aim to show the identity for $G^{(II) }_{uu,j},\ j=1,2.$

From \eqref{hatGv} , we now have
\begin{equation} \label{Guu2}
\int_{\Omega} (G^{(II)}_{uu,1}-G^{(II)}_{uu,2}) h^2(x;\lambda)  e^{-\frac{i}{\sqrt{d_2}}\xi\cdot x} dx=0. 
\end{equation}

Since $h(x;\lambda)$ can be chosen to be a non-zero eigenfunction of $\Delta,$ hence
\[G^{(II) }_{uu,1}=G^{(II) }_{uu,2}.\]

\textbf{\textit{Recovery of $G^{(II) }_{vv,j}$} } Similarly, we can recover $G^{(II) }_{vv,j}$ given  $G^{(II) }_{vu},G^{(II) }_{uv},$ and $G^{(II) }_{uu}.$

Substituting the given information into \eqref{hatGv} and we have, 
\begin{equation} \label{Gvv2}
\int_{Q} (G^{(II)}_{vv,1}-G^{(II)}_{vv,2})l^2(x;\mu)e^{-\frac{i}{\sqrt{d_2}}\xi\cdot x} dx=0.  
\end{equation}

Since $l(x;\mu)$ can be chosen to be a non-zero eigenfunction of $\Delta,$ hence
\[G^{(II) }_{vv,1}=G^{(II) }_{vv,2}.\]

Finally, we can use mathematical induction and repeat similar arguments in the above two subsections to show that for each $k,$ $k \geq 3, k \in \mathbb{N} ,$ we can first show 
\[F^{(k)}_{1}-F^{(k)}_{2}=0,\]
where $F^{(k)}_{j}, j=1,2$ covers all the partial derivatives related to $u$ and $v,$ for example $F^{(k)}_{uu\cdots u,j},F^{(k)}_{vu\cdots u,j},F^{(k)}_{vv\cdots u,j}, \dots$ and so on.

Then we can simultaneously recover 
\[G^{(k)}_{1}-G^{(k)}_{2}=0,\]
for each $k,$ $k \geq 3, k \in \mathbb{N} ,$ where $G^{(k)}_{j}, j=1,2$ covers all the partial derivatives related to $u$ and $v.$
Hence, it reflects that
\[F_1(x,t,u,v)=F_2(x,t,u,v), \  G_1(x,t,u,v)=G_2(x,t,u,v)\  \text{in}\  \Omega \times \mathbb{R}. \]

The proof is complete.$ \hfill{\square}$

Through the consideration of successively higher orders of linearization, we are able to successively recover the Taylor coefficients of $F$ and $G$ simultaneously.

\begin{rmk}\label{Differecestate}
\textit{ 
In \cite{LL23}, the analytic functions}  $F(x,t,p,q), G(x,t,p,q):\Omega\times (0,T)\times \mathbb{R}\times \mathbb{R} \to \mathbb{R}$ \textit{are of the form}  
\[F(x,t,p,q):=\sum\limits^{\infty}_{m+n\geq3}\alpha_{mn}(x,t) p^m q^n,\  G(x,t,p,q):=\sum\limits^{\infty}_{m+n\geq1}\beta_{mn}(x,t) p^m q^n,\]
\textit{where} $m,n\geq 0.$ \textit{Under several strict restrictions, they can only recover the} $\alpha_{m0}$  \textit{and}  $\beta_{m0}$ \textit{for} $m \geq 2.$ \textit{For \cite{LZ23}, the model considered by the authors involved source functions of the type
} $F(x,v)$ depending only on the second variable and is independent of time.
\textit{Our recovery items, on the other hand, are subject to significantly looser assumptions in comparison and can adapt to more complex forms for}  $F(x,t,p,q)$ \textit{and} $ G(x,t,p,q).$ \textit{We can recover not only} $F_{u\cdots u}, G_{u\cdots u}$ \textit{terms, but also any partial derivatives in a higher-order variation form, which provides much wider options in applications compared to \cite{LL23} and \cite{LZ23}.}
\end{rmk}

{\centering \section{Biological Applications:  diffusive Bazykin model}  \label{appbio}}

Our results can be applied to a variety of models, especially in biology. In this section, we explain this by describing how we simultaneously recover the coefficients in \eqref{hydraeq} and \eqref{Difex}. 

\subsection{Application for hydra effect}
First, we consider \eqref{hydraeq} in the following form:
\begin{equation} \label{apbiohydra}
\begin{cases} 
\partial_{t} u-d_1\Delta u= u(a-b-eu)-\tilde{F}(x,t,u,v) &\  \text{in}\  Q,\\ 
\partial_{t} v-d_2\Delta v= -m v+\mu \lambda v^2+\tilde{G}(x,t,u,v) &\  \text{in}\  Q,\\ 
\partial_{\nu} u=\partial_{\nu} v=0 &\  \text{on}\  \Sigma,\\ 
u(x,0)= f(x) , \; v(x,0)= g(x)& \  \text{in} \ \Omega. 
\end{cases} 
\end{equation} 
where $\tilde{F}(x,t,u,v)=(p+\lambda v)uv,\tilde{G}_(x,t,u,v)=\mu p uv.$ We recall that this is the model for hydra effect, where the well-posedness of \eqref{apbiohydra} is discussed in \cite{LBi22} and Section \ref{forp}. It is apparent to see that $(0,0)$ is a constant solution to \eqref{apbiohydra}, and $\tilde{F}(x,t,u,v), \tilde{G}(x,t,u,v)$ are analytic with respect to $u$ and $v.$ We can give admissible classes for $\tilde{F},\tilde{G}$ as 

\begin{defi}\label{admHy1}
  We say that $U(x,t,p,q): \mathbb{R}^n \times \mathbb{R} \times \mathbb{C}  \times \mathbb{C} \to \mathbb{C}$ is admissible, denoted by $U \in \mathcal{C},$ if:

(a) The map $(p,q) \mapsto U(\cdot, \cdot, p,q) $ is holomorphic with value in $C^{2+\alpha, 1+\frac{\alpha}{2}}(\bar{Q}),$

(b) $U(x,t,0,0)=0$ for all $(x,t) \in Q,$

(c) $U^{(0,1)}(x,t,\cdot,0)=U^{(1,0)}(x,t,0,\cdot)=0$ for all $(x,t) \in Q,$ 

(d) Taylor coefficients of all orders for $U$ are constants.  
\end{defi}

\begin{defi}\label{admHy2}
    We say that $V(x,t,p,q): \mathbb{R}^n \times \mathbb{R} \times \mathbb{C}  \times \mathbb{C} \to \mathbb{C}$ is admissible, denoted by $V \in \mathcal{D},$ if:

(a) The map $(p,q) \mapsto V(\cdot, \cdot, p,q) $ is holomorphic with value in $C^{2+\alpha, 1+\frac{\alpha}{2}}(\bar{Q}),$

(b)$V(x,t,0,0)=0$ for all $(x,t) \in Q,$ 

(c) $V^{(1,0)}(x,t,0,\cdot)=0$ for all $(x,t) \in Q,$ 

(d) Taylor coefficients of all orders for $V$ are constants.
\end{defi}

The measurement map applied in this system is 
\[\mathcal{M}^{+}_{F,G}(f(x),g(x))=\big( (u(x,t),v(x,t))\vert_{\Sigma},u(x,T),v(x,T)\big),\    x \in \Omega.\]
Our inverse problem is to recover the unknown coefficients of $\tilde{F}$ and $\tilde{G}$ given two identical measurement maps. Then we can obtain the following result.

\begin{prop}\label{appprop1}
Assume $\tilde{F}_j \in \mathcal{C},$ $\tilde{G}_j \in \mathcal{D},$  $ (j=1,2). $ Let $\mathcal{M}^{+}_{\tilde{F}_j, \tilde{G}_j}$ be the measurement map associated to the following system:
\begin{equation} \label{apbiohydra2}
\begin{cases} 
\partial_{t} u_{j}-d_1\Delta u_j= u_j(a-b-eu_j)-\tilde{F}_j(x,t,u,v) &\  \text{in}\  Q,\\ 
\partial_{t} v_{j}-d_2\Delta v_j= -m v_j+\mu \lambda v^2_j+\tilde{G}_j(x,t,u,v) &\  \text{in}\  Q,\\ 
\partial_{\nu} u_{j}=\partial_{\nu} v_{j}=0 &\  \text{on}\  \Sigma,\\ 
u_j(x,0)= f(x) , \; v_j (x,0)= g(x)& \  \text{in} \ \Omega.
\end{cases} 
\end{equation} 
If for any $f,g \in C^{2+\alpha}(\Omega),$ one has  \[\mathcal{M}^{+}_{\tilde{F}_1,\tilde{G}_1}(f,g)=\mathcal{M}^{+}_{\tilde{F}_2,\tilde{G}_2}(f,g),\]
then it holds that 
\[ \tilde{F}_1(x,t,u,v)=\tilde{F}_2(x,t,u,v),\  \, \tilde{G}_1(x,t,u,v)=\tilde{G}_2(x,t,u,v)\  \ \text{in}\  \, \Omega \times \mathbb{R}.\]
\end{prop}

\textit{Proof. } 
We use Theorem \ref{mainthm} directly for the proof. 

Since $(0,0)$ is a solution of \eqref{apbiohydra}, we fix $u_0=0, v_0=0.$ It is apparent to see 
\[F_j(x,t,u,v)=u_j(a-b-eu_j)-\tilde{F}_j(x,t,u,v) \in \mathcal{A} ,\]
\[G_j(x,t,u,v)=-m v_j+\mu \lambda v^2_j+\tilde{G}_j(x,t,u,v) \in \mathcal{B} ,\]
satisfy the assumptions in Theorem \ref{mainthm}. 

Therefore, we have the result for Proposition \ref{appprop1}. $ \hfill{\square}$

Physically, our result can be interpreted in the following sense. We input the population densities of both prey and predator at the initial time of a field trip or an experiment and measure their population density values at the boundary of a bounded region $\Omega$ as well as at a given end time. If all these values are the same, there is only one possible value for the attack rate of an individual predator, the efficiency of food conversion into offspring, and the strength of cooperation during hunting, respectively, which we could recover sequentially.

\subsection{Application for Holling-Tanner type}
The next model we consider is the Holling-Tanner type system, which we recall models describing a specific prey preference for one predator. It is given by the following:
\begin{equation} \label{Difuse}
\begin{cases} 
\partial_{t}u-d_1\Delta u= u(1-u)-\hat{F}(x,t,u,v) &\  \text{in} \   Q,\\ 
\partial_{t}v-d_2\Delta v= -\frac{c}{a} v-v^2+\hat{G}(x,t,u,v) &\  \text{in}  \  Q,\\ 
\partial_{\nu}u=\partial_{\nu}v=0&\  \text{on}\   \Sigma,\\ 
u(x,0)= f(x) , \; v (x,0)= g(x)&\  \text{in}\   \Omega. 
\end{cases} 
\end{equation} 
In the following, we denote $\frac{c}{a}$ as $\delta.$ The well-posedness of \eqref{Difuse} is discussed in \cite{LBi22} and Section \ref{forp}. It is trivial to see that $(1,0)$ and $(0,0)$ are two solutions to \eqref{Difuse}. The process of simultaneous recovery is similar for the chosen initial inputs $(1,0)$ and $(0,0),$ so we only illustrate the case $(1,0)$ as an example. Choosing $(1,0)$ as the initial data can also display the advantage of applying high-order variation method.

Known from \eqref{Difex}, 
\[\hat{F}(x,t,u,v)=\frac{\beta uv}{\alpha+u},\  \hat{G}(x,t,u,v)=\frac{\gamma uv}{\alpha+u}\  \text{with}\  \alpha=\frac{1}{AK},\, \beta=\frac{b}{AhK}, \, \gamma=\frac{d}{Aa} \]
can be easily seen to be analytic with respect to $u$ and $v,$ by using the Taylor series expansion. Therefore we can define the admissible class directly for $\hat{F}$ and $\hat{G}.$ 

\begin{defi}\label{admHT1}
    We say that $U(x,t,p,q): \mathbb{R}^n \times \mathbb{R} \times \mathbb{C}  \times \mathbb{C} \to \mathbb{C}$ is admissible, denoted by $U \in \mathcal{E},$ if:

(a) The map $(p,q) \mapsto U(\cdot, \cdot, p,q) $ is holomorphic with value in $C^{2+\alpha, 1+\frac{\alpha}{2}}(\bar{Q}),$

(b) $U(x,t,1,0)=0$ for all $(x,t) \in Q,$

(c) $U^{(0,1)}(x,t,\cdot,0)=U^{(1,0)}(x,t,1,\cdot)=0$ for all $(x,t) \in Q,$ 

(d) Taylor coefficients of all orders for $U$ are constants.
\end{defi}

\begin{defi}\label{admHT2}
    We say that $V(x,t,p,q): \mathbb{R}^n \times \mathbb{R} \times \mathbb{C}  \times \mathbb{C} \to \mathbb{C}$ is admissible, denoted by $V \in \mathcal{F},$ if:

(a) The map $(p,q) \mapsto V(\cdot, \cdot, p,q) $ is holomorphic with value in $C^{2+\alpha, 1+\frac{\alpha}{2}}(\bar{Q}),$

(b)$V(x,t,1,0)=0$ for all $(x,t) \in Q,$ 

(c) $V^{(1,0)}(x,t,1,\cdot)=0$ for all $(x,t) \in Q,$ 

(d) Taylor coefficients of all orders for $V$ are constants.
\end{defi}

The measurement map we use is 
\[\mathcal{M}^{+}_{\hat{F},\hat{G}}(f(x),g(x))=\big( (u(x,t),v(x,t))\vert_{\Sigma},u(x,T),v(x,T)\big),\    x \in \Omega. \]
Our inverse problem is to recover the unknown coefficients of $\hat{F}$ and $\hat{G}$ given two identical measurement maps. Then we can obtain the following result.

\begin{prop}\label{appHT}
     Assume $\hat{F}_j \in \mathcal{E},$ $\hat{G}_j \in \mathcal{F},$  $ (j=1,2). $ Let $\mathcal{M}^{+}_{\hat{F}_j, \hat{G}_j}$ be the measurement map associated to the following system:
\begin{equation} \label{apbio}
\begin{cases} 
\partial_{t} u_{j}-d_1\Delta u_j= u_j(1-u_j)-\hat{F}_j(x,t,u,v) &\  \text{in}\  Q,\\ 
\partial_{t} v_{j}-d_2\Delta v_j= -\delta v_j-v^2_j+\hat{G}_j(x,t,u,v) &\  \text{in}\  Q,\\ 
\partial_{\nu} u_{j}=\partial_{\nu} v_{j}=0 &\  \text{on}\  \Sigma,\\ 
u_j(x,0)= f(x) , \; v_j (x,0)= g(x)& \  \text{in} \ \Omega.
\end{cases} 
\end{equation} 
If for any $f,g \in C^{2+\alpha}(\Omega),$ one has  \[\mathcal{M}^{+}_{\hat{F}_1,\hat{G}_1}(f,g)=\mathcal{M}^{+}_{\hat{F}_2,\hat{G}_2}(f,g),\]
then it holds that 
\[ \hat{F}_1(x,t,u,v)=\hat{F}_2(x,t,u,v),\  \, \hat{G}_1(x,t,u,v)=\hat{G}_2(x,t,u,v)\  \ \text{in}\  \, \Omega \times \mathbb{R}.\]
\end{prop}

\textit{Proof. } 
We use Theorem \ref{mainthm} directly for the proof.

Since  $(1,0)$ is a solution of \eqref{apbio}, we fix $u_0=1, v_0=0.$ It is apparent to see 
\[F_j(x,t,u,v)=ku_j(1-u_j)-\hat{F}_j(x,t,u,v) \in \mathcal{A} ,\]
\[G_j(x,t,u,v)=-\delta v_j-v^2_j+\hat{G}_j(x,t,u,v) \in \mathcal{B} ,\]
satisfy the assumptions in Theorem \ref{mainthm}. 

Therefore, we have the result for Proposition \ref{appHT}. $ \hfill{\square}$

Once again, we can interpret our result in physical meaning as below. We input the population densities of prey and predator at the initial time of a field trip or an experiment and measure their population density values at the boundary of a bounded region $\Omega$ and at a given end time. Suppose all these values are the same, and the maximum consumption of the predator as well as the carrying capacity are known. In that case, there is only one possible value for the coefficient of predator competition for the crowding effect and the conversion efficiency from the prey to the predator. Alternatively, by knowing the last two terms of the mentioned coefficients, we can obtain the only possibility value for the maximum consumption of the predator and the carrying capacity.

\subsection{Application for the classic L-V model}
In particular case where the system \eqref{Difuse} is of the form \eqref{classic} corresponding to the Bazykin model, the result applies, and is given as follows.

\begin{prop}\label{appgeneral}
    Assume $F_j \in \mathcal{A},$ $G_j \in \mathcal{B},$  $ (j=1,2). $ Let $\mathcal{M}^{+}_{F_j, G_j}$ be the measurement map associated to the following system:
\begin{equation} \label{apbio2}
\begin{cases} 
\partial_{t} u_{j}=au_j(1-\frac{u_j}{K})-v_jp(u_j,v_j) &\  \text{in}\  Q,\\ 
\partial_{t} v_{j}=v_j(-c+dp(u_j,v_j))-hv^2_j &\  \text{in}\  Q,\\ 
\partial_{\nu} u_{j}=\partial_{\nu} v_{j}=0 &\  \text{on}\  \Sigma,\\ 
u_j(x,0)= f(x) , \; v_j (x,0)= g(x)& \  \text{in} \ \Omega, 
\end{cases} 
\end{equation} 
where the functions $F$ and $G$ are specifically given by
\[F_j(x,t,u,v)=au_j(1-\frac{u_j}{K})-v_jp(u_j,v_j),\]
and
\[G_j(x,t,u,v)=v_j(-c+dp(u_j,v_j))-hv^2_j.\]
If for any $f,g \in C^{2+\alpha}(\Omega),$ one has  \[\mathcal{M}^{+}_{F_1,G_1}(f,g)=\mathcal{M}^{+}_{F_2,G_2}(f,g),\]
then it holds that 
\begin{equation}\label{biocon}
    F_1(x,t,u,v)=F_2(x,t,u,v),\  \, G_1(x,t,u,v)=G_2(x,t,u,v)\  \ \text{in}\  \, \Omega \times \mathbb{R}.
\notag\end{equation}
\end{prop}

We can verify Proposition \ref{appgeneral} directly from Theorem \ref{mainthm}. The proposition tells us that if all the inputs population densities of prey and predator at the initial time of a field trip or an experiment and all the measurements of their population density values at the boundary as well as at the given end time are identical, we can obtain the only probability for the coefficient depicting the predator competition caused by self-limitation and the Taylor coefficients of the particular functional response.

Moreover, instead of setting a non-negative initial predator population density, we can assume a strictly positive initial data $g_1(x,t)>0.$ In its physical meaning, we know there is a certain population of predator in the observation area, which sounds more reasonable for a biological experiment.

Note that apart from a broader range for initial data chosen compared to \cite{LL23}, there is a difference in background settings. Though both of us measure the flux at prey and predator population density at the boundary of a bounded encircled region $\Omega,$   \cite{LL23} inputs a particular population of prey from the boundary into $\Omega,$ while this paper uses the population densities at an arbitrary cut-in time. Our model offers a more versatile method for observers and experimenters in a project involving a macro-scale ecological background. 

\medskip

	\noindent\textbf{Acknowledgment.} 
	The work was supported by the Hong Kong RGC General Research Funds (projects 11311122, 11300821 and 12301420),  the NSFC/RGC Joint Research Fund (project N\_CityU101/21), and the ANR/RGC Joint Research Grant, A\_CityU203/19.

\bibliographystyle{plain}
\bibliography{newre}

\begin{thebibliography}{10}

\bibitem{AbrHy09}
Peter~A Abrams.
\newblock When does greater mortality increase population size{?} {The} long history and diverse mechanisms underlying the hydra effect.
\newblock {\em Ecology letters}, 12(5):462--474, 2009.

\bibitem{AbrHy05}
Peter~A Abrams and Hiroyuki Matsuda.
\newblock The effect of adaptive change in the prey on the dynamics of an exploited predator population.
\newblock {\em Canadian Journal of Fisheries and Aquatic Sciences}, 62(4):758--766, 2005.

\bibitem{B76}
Alexander~D Bazykin.
\newblock Structural and dynamic stability of model predator-prey systems.
\newblock 1976.

\bibitem{B98}
Alexander~D Bazykin.
\newblock {\em Nonlinear dynamics of interacting populations}.
\newblock World Scientific, 1998.

\bibitem{BAC09}
Assia Benabdallah, Michel Cristofol, Patricia Gaitan, and Masahiro Yamamoto.
\newblock Inverse problem for a parabolic system with two components by measurements of one component.
\newblock {\em Applicable Analysis}, 88(5):683--709, 2009.

\bibitem{CSE22}
Salah-Eddine Chorfi and Lahcen Maniar.
\newblock Stable determination of coefficients in semilinear parabolic system with dynamic boundary conditions.
\newblock {\em Inverse Problems}, 38(11):115007, 2022.

\bibitem{DL23mfg}
Ming-Hui Ding, Hongyu Liu, and Guang-Hui Zheng.
\newblock Determining a stationary mean field game system from full/partial boundary measurement.
\newblock {\em arXiv preprint arXiv:2308.06688}, 2023.

\bibitem{DL23}
Ming-Hui Ding, Hongyu Liu, and Guang-Hui Zheng.
\newblock On inverse problems for several coupled {PDE} systems arising in mathematical biology.
\newblock {\em Journal of Mathematical Biology}, 2023.

\bibitem{Du97}
Yihong Du and Yuan Lou.
\newblock Some uniqueness and exact multiplicity results for a predator-prey model.
\newblock {\em Transactions of the American Mathematical Society}, 349(6):2443--2475, 1997.

\bibitem{DY1}
Yihong Du and Yuan Lou.
\newblock Qualitative behaviour of positive solutions of a predator-prey model: effects of saturation.
\newblock {\em Proceedings of the Royal Society of Edinburgh Section A: Mathematics}, 131(2):321--349, 2001.

\bibitem{EW9}
Heinz~W Engl, Christoph Flamm, Philipp K{\"u}gler, James Lu, Stefan M{\"u}ller, and Peter Schuster.
\newblock Inverse problems in systems biology.
\newblock {\em Inverse Problems}, 25(12):123014, 2009.

\bibitem{F97}
Herbert~I Freedman.
\newblock Stability analysis of a predator-prey system with mutual interference and density-dependent death rates.
\newblock {\em Bulletin of Mathematical Biology}, 41:67--78, 1979.

\bibitem{F980ec}
Herbert~I Freedman.
\newblock Deterministic mathematical models in population ecology.
\newblock {\em (No Title)}, 1980.

\bibitem{FA8}
Avner Friedman.
\newblock {\em Partial differential equations of parabolic type}.
\newblock Courier Dover Publications, 2008.

\bibitem{Ha92}
Josef Hainzl.
\newblock Multiparameter bifurcation of a predator-prey system.
\newblock {\em SIAM journal on mathematical analysis}, 23(1):150--180, 1992.

\bibitem{H09ratio}
Mainul Haque.
\newblock Ratio-dependent predator-prey models of interacting populations.
\newblock {\em Bulletin of mathematical biology}, 71:430--452, 2009.

\bibitem{Haq11BD}
Mainul Haque.
\newblock A detailed study of the {B}eddington--{D}eangelis predator--prey model.
\newblock {\em Mathematical Biosciences}, 234(1):1--16, 2011.

\bibitem{J21glo}
Xin Jiang, Zhikun She, and Shigui Ruan.
\newblock Global dynamics of a predator-prey system with density-dependent mortality and ratio-dependent functional response.
\newblock {\em Discrete Contin. Dyn. Syst. Ser. B}, 26(4):1967--1990, 2021.

\bibitem{Ku98}
Yuri~A Kuznetsov, Iu~A Kuznetsov, and Y~Kuznetsov.
\newblock {\em Elements of applied bifurcation theory}, volume 112.
\newblock Springer, 1998.

\bibitem{LSU88}
Olga~Aleksandrovna Ladyzhenskaia, Vsevolod~Alekseevich Solonnikov, and Nina~N Ural'tseva.
\newblock {\em Linear and quasi-linear equations of parabolic type}, volume~23.
\newblock American Mathematical Soc., 1988.

\bibitem{L78}
Anthony Leung.
\newblock Limiting behaviour for a prey-predator model with diffusion and crowding effects.
\newblock {\em J. Math. Biol}, 6(1):87--93, 1978.

\bibitem{SR}
Yi-Hsuan Lin, Hongyu Liu, Xu~Liu, and Shen Zhang.
\newblock Simultaneous recoveries for semilinear parabolic systems.
\newblock {\em Inverse Problems}, 38(11):115006, 2022.

\bibitem{LL23}
Hongyu Liu and Catharine~WK Lo.
\newblock Determining a parabolic system by boundary observation of its non-negative solutions with applications.
\newblock {\em Inverse problems}, 2023.

\bibitem{LMZ23}
Hongyu Liu, Chenchen Mou, and Shen Zhang.
\newblock Inverse problems for mean field games.
\newblock {\em Inverse Problems, Volume 39, Number 8}, 2023.

\bibitem{LZ22}
Hongyu Liu and Shen Zhang.
\newblock On an inverse boundary problem for mean field games.
\newblock {\em arXiv preprint arXiv:2212.09110}, 2022.

\bibitem{LZ23}
Hongyu Liu and Shen Zhang.
\newblock Simultaneously recovering running cost and {H}amiltonian in mean field games system.
\newblock {\em arXiv preprint arXiv:2303.13096}, 2023.

\bibitem{L1925}
Alfred~James Lotka.
\newblock {\em Elements of physical biology}.
\newblock Williams \& Wilkins, 1925.

\bibitem{LGl21}
Min Lu and Jicai Huang.
\newblock Global analysis in {B}azykin's model with {H}olling {II} functional response and predator competition.
\newblock {\em Journal of Differential Equations}, 280:99--138, 2021.

\bibitem{LBi22}
Min Lu, Chuang Xiang, Jicai Huang, and Hao Wang.
\newblock Bifurcations in the diffusive {B}azykin model.
\newblock {\em Journal of Differential Equations}, 323:280--311, 2022.

\bibitem{PV}
Chia-Ven Pao.
\newblock {\em Nonlinear parabolic and elliptic equations}.
\newblock Springer Science \& Business Media, 2012.

\bibitem{PW83}
Stephen~W Provencher and Robert~H Vogel.
\newblock Regularization techniques for inverse problems in molecular biology.
\newblock In {\em Numerical Treatment of Inverse Problems in Differential and Integral Equations: Proceedings of an International Workshop, Heidelberg, Fed. Rep. of Germany, August 30—September 3, 1982}, pages 304--319. Springer, 1983.

\bibitem{Reva}
Max Rietkerk, Robbin Bastiaansen, Swarnendu Banerjee, Johan van~de Koppel, Mara Baudena, and Arjen Doelman.
\newblock Evasion of tipping in complex systems through spatial pattern formation.
\newblock {\em Science}, 374(6564):eabj0359, 2021.

\bibitem{SL20}
Danxia Song, Chao Li, and Yongli Song.
\newblock Stability and cross-diffusion-driven instability in a diffusive predator-prey system with hunting cooperation functional response.
\newblock {\em Nonlinear Analysis: Real World Applications}, 54:103106, 2020.

\bibitem{T52}
Alan Turing.
\newblock The chemical basis of morphogenesis.
\newblock {\em B Jack Copeland}, 519, 1952.

\bibitem{V1926}
Vito Volterra.
\newblock Fluctuations in the abundance of a species considered mathematically.
\newblock {\em Nature}, 118(2972):558--560, 1926.

\end{thebibliography}

\end{document}